\documentclass[reqno,11pt]{amsart} 
\usepackage{amsmath}
\usepackage{amssymb}
\usepackage{amsthm}
\usepackage{verbatim}
\usepackage{xcolor}
\usepackage{graphics}

\newcommand\NoBlackBoxes{\global\overfullrule0pt}
\NoBlackBoxes
\theoremstyle{plain} 

\def\4{\kern1pt}

\def\6{\vphantom0}

\def\8{\kern-10pt}
\def\7#1{_{(#1)}}

\begin{document}

\title{ASYMPTOTIC EXPANSIONS AND TWO-SIDED BOUNDS \\
IN RANDOMIZED CENTRAL LIMIT THEOREMS}

\author{S. G. Bobkov$^{1,3}$}
\thanks{1) 
School of Mathematics, University of Minnesota, USA}

\author{G. P. Chistyakov$^{2,3}$}

\author{F. G\"otze$^{2,3}$}
\thanks{2) 
Faculty of Mathematics, University of Bielefeld, Germany}
\thanks{3) Research supported by SFB 1283, BSF grant 2016050, 
and NSF grant DMS-2154001}

\dedicatory{Dedicated to the memory of Gennadiy P. Chistyakov\;  *\,May 1, 1945 \; $\dagger$\,December 30, 2022.}

\subjclass
{Primary 60E, 60F} 
\keywords{Typical distributions, normal approximation,
central limit theorem} 

\begin{abstract}
Lower and upper bounds are explored for 
the uniform (Kolmogorov) and $L^2$-distances between the distributions 
of weighted sums of dependent summands and the normal law. 
The results are illustrated for several classes of random variables 
whose joint distributions are supported on Euclidean spheres.
We also survey several results on improved rates of normal approximation
in randomized central limit theorems.
\end{abstract}

\maketitle
\markboth{S. G. Bobkov, G. P. Chistyakov and F. G\"otze}{
CLT for weighted sums}

\def\theequation{\thesection.\arabic{equation}}
\def\E{{\mathbb E}}
\def\R{{\mathbb R}}
\def\C{{\mathbb C}}
\def\P{{\mathbb P}}
\def\I{{\mathbb I}}
\def\S{{\mathbb S}}

\def\s{{\mathfrak s}}

\def\H{{\rm H}}
\def\Im{{\rm Im}}
\def\Tr{{\rm Tr}}

\def\k{{\kappa}}
\def\M{{\cal M}}
\def\Var{{\rm Var}}
\def\Ent{{\rm Ent}}
\def\O{{\rm Osc}_\mu}

\def\ep{\varepsilon}
\def\phi{\varphi}
\def\F{{\cal F}}
\def\L{{\cal L}}

\def\be{\begin{equation}}
\def\en{\end{equation}}
\def\bee{\begin{eqnarray*}}
\def\ene{\end{eqnarray*}}


\section{{\bf Introduction}}
\setcounter{equation}{0}

\vskip2mm
\noindent
A random vector $X = (X_1,\dots,X_n)$ in $\R^n$ ($n \geq 2$) defined
on the probability space $(\Omega,{\mathfrak F},\P)$ is called isotropic, if 
$$
\E X_i X_j = \delta_{ij} \quad {\rm for \ all} \ i,j \leq n,
$$ where $\delta_{ij}$ 
is the Kronecker symbol. Equivalently, all weighted sums
$$
S_\theta = \theta_1 X_1 + \dots + \theta_n X_n, \qquad 
\theta = (\theta_1,\dots,\theta_n), \ \ \theta_1^2 + \dots + \theta_n^2=1,
$$
with coefficients from the unit sphere $\S^{n-1}$ in $\R^n$ have a second
moment $\E S_\theta^2 = 1$. In this case, provided that the Euclidean norm
$|X|$ is almost constant, and if $n$ is large, a theorem due 
to Sudakov \cite{Su} asserts that the distribution functions 
$$
F_\theta(x) = \P\{S_\theta \leq x\}, \quad x \in \R,
$$
are well approximated for most of $\theta \in \S^{n-1}$ by the standard 
normal distribution function
$$
\Phi(x) = \frac{1}{\sqrt{2\pi}} \int_{-\infty}^x e^{-y^2/2}\,dy.
$$
Here, ``most" should refer to the normalized Lebesgue measure $\s_{n-1}$ 
on the sphere. This property may be quantified, for example, in terms 
of the Kolmogorov distance
$$
\rho(F_\theta,\Phi) = \sup_x |F_\theta(x) - \Phi(x)|.
$$

Being rather universal (since no independence of the components $X_k$ is 
required), randomized central limit theorems of such type have received 
considerable interest in recent years. For the history, bibliography, and 
interesting connections with other concentration problems we refer 
an interested reader to \cite{B-C-G3}, \cite{B-C-G4}, \cite{B-G-V-V}. 
Let us mention one general upper bound
\be
\E_\theta\, \rho(F_\theta,\Phi) \, \leq \, 
c\,(1 + \sigma_4)\,\frac{\log n}{\sqrt{n}},
\en
which holds true with an absolute constant $c>0$ for any isotropic random 
vector $X$ (cf. Theorem 1.2 in \cite{B-C-G3}). Here and elsewhere, 
$\E_\theta$ denotes an integral over $\S^{n-1}$ with respect to 
the measure $\s_{n-1}$, and the bound involves the variance-type functional 
$$
\sigma_4^2 = \sigma_4^2(X) = \frac{1}{n}\,\Var(|X|^2) \quad (\sigma_4 \geq 0).
$$ 

Modulo a logarithmic factor, the bound (1.1) exhibits a standard rate 
of normal approximation for $F_\theta$, in analogy with the classical case 
of independent identically distributed (iid) summands with equal coefficients. 
It turns out, however, that in the model with arbitrary $\theta \in \S^{n-1}$ 
and independent components $X_k$, the standard rate for $\rho(F_\theta,\Phi)$ 
is dramatically improved to the order $1/n$ on average and actually for most 
of $\theta$. Motivated by the seminal paper of Klartag and Sodin \cite{K-S},
this interesting phenomenon was recently studied in \cite{B-C-G4}, \cite{B-C-G5} 
for dependent data under certain correlation-type conditions. The last chapters
of this paper provide a short account of these improved
rates of normal approximation.

One of the main aims of this work is to develop lower bounds with a similar 
standard rate as in (1.1) (modulo logarithmic factors) and to illustrate them
with a number of examples of random variables $X_k$ often appearing
in Functional Analysis. These results 
rely on a careful examination of the closely related $L^2$-distance
$$
\omega(F_\theta,\Phi) = 
\Big(\int_{-\infty}^\infty (F_\theta(x) - \Phi(x))^2\,dx\Big)^{1/2}.
$$
Similarly to (1.1), it can be shown that for the class of isotropic 
random vectors the inequality
\be
\E_\theta\, \omega^2(F_\theta,\Phi) \, \leq \, c\,(1 + \sigma_4^2)\,\frac{1}{n}
\en
holds without an unnecessary logarithmic term. However, in order to 
explore the real behavior of the average $L^2$-distance, some other 
characteristics of the distribution of $X$ are required. For example, 
assuming that the distribution is supported on the sphere 
$\sqrt{n}\,\S^{n-1}$, the $L^2$-distance admits an asymptotic expansion 
in terms of the moment functionals (normalized $L^p$-norms)
$$
m_p = m_p(X) = \frac{1}{\sqrt{n}}\,\big(\E \left<X,Y\right>^p\big)^{1/p}
= \frac{1}{\sqrt{n}}\,
\Big(\sum\, (\E X_{i_1} \dots X_{i_p})^2\Big)^{1/p}.
$$
Here, $Y$ is an independent copy of $X$, and the summation is performed
over all indices $1 \leq i_1,\dots,i_p \leq n$. The second representation
shows that these functionals are non-negative for any integer $p \geq 1$. 
Note that $m_1 = 0$ if $X$ has mean zero, $m_2 = 1$ if $X$
is isotropic, and $m_p = 0$ with odd $p$ when
the distribution of $X$ is symmetric about the origin.
The following expansion involves the moments $m_p$ up to order 4.

\vskip5mm
{\bf Theorem 1.1.} {\sl Let $X$ be an isotropic random vector 
in $\R^n$ with mean zero and such that $|X|^2 = n$ a.s. We have
\be
\E_\theta\, \omega^2(F_\theta,\Phi)  \, = \,
\frac{c}{n^{3/2}}\, m_3^3 + O\Big(\frac{1}{n^2}\,m_4^4\Big)
\en
with $c = \frac{1}{16 \sqrt{\pi}}$.  Similarly, with some absolute 
constants $c_1,c_2>0$,
\be
\E_\theta\, \rho^2(F_\theta,\Phi)  \, \leq \,
\frac{c_1 \log n}{n^{3/2}}\, m_3^3 + 
\frac{c_2 (\log n)^2}{n^2}\, m_4^4.
\en
}

\vskip2mm
As we will see, in the general isotropic case without the support 
assumption, but with bounded $\sigma_4$, the average $L^2$-distance 
is described by a more complicated formula
\begin{eqnarray}
\E_\theta\, \omega^2(F_\theta,\Phi) 
 & = &
\frac{1}{\sqrt{2\pi n}} \, \Big(1 + \frac{1}{8n}\Big)\, \E\sqrt{|X|^2 + |Y|^2}
\nonumber \\
 & & - \ \frac{1}{\sqrt{2\pi n}} \, 
\Big(1 + \frac{1}{4n}\Big)\,\E\,|X-Y| + 
O\Big(\frac{1 + \sigma_4^2}{n^2}\Big),
\end{eqnarray}
which holds whenever $\E\,|X|^2 = n$.

In the setting of Theorem 1.1, using the pointwise bound 
$|\left<X,Y\right>| \leq n$ together with the isotropy assumption, 
we have $\E \left<X,Y\right>^3 \leq n^2$ and 
$\E \left<X,Y\right>^4 \leq n^3$. Therefore, the inequalities 
(1.3)-(1.4) yield with some absolute constant $c>0$
\be
\E_\theta\, \omega^2(F_\theta,\Phi) \leq \frac{c}{n}, \quad
\E_\theta\, \rho^2(F_\theta,\Phi) \leq \frac{c\, (\log n)^2}{n},
\en
thus recovering the upper bounds (1.1)-(1.2) for this particular case 
(since $\sigma_4 = 0$). On the other hand, for a large variety of 
examples, such bounds turn out to be optimal and may be reversed 
modulo a logarithmic factor (for large~$n$).
To see this, one may use the following lower bound which will be
derived from a slightly modified variant of (1.5).

\vskip5mm
{\bf Theorem 1.2.} {\sl Let $X$ be a random vector in $\R^n$ 
satisfying $\E\,|X|^2 = n$, and let $Y$ be its independent copy. 
For some absolute constants $c_1,c_2 > 0$, we have
\be
\E_\theta\, \omega^2(F_\theta,\Phi) \ \geq \ 
c_1\,\P\Big\{|X-Y| \leq \frac{1}{2}\sqrt{n}\,\Big\} - 
c_2 \frac{1 + \sigma_4^4}{n^2}.
\en
}

\vskip2mm
Thus, if the probability in (1.7) is of order at least $1/n$, and 
$\sigma_4$ is bounded, the right-hand side of this bound will be of 
the same order. If, for example, $|X| = \sqrt{n}$ a.s., we then obtain 
that $\E_\theta\, \omega^2(F_\theta,\Phi) \sim 1/n$.
In order to derive a similar conclusion for the Kolmogorov distance, 
one may refer to the next statement.

\vskip5mm
{\bf Theorem 1.3.} {\sl Let $X$ be an isotropic random vector in 
$\R^n$ such that $|X| \leq b\sqrt{n}$ a.s. Suppose that we have a lower 
bound at the standard rate
$$
\E_\theta\,\omega^2(F_\theta,\Phi) \, \geq \, 
\frac{D}{n}
$$
with some $D>0$. Then with some absolute constants $c_0,c_1 > 0$
$$
\E_\theta\,\rho(F_\theta,F)\,\geq\,  
\frac{c_0}{(1 + \sigma_4)^3\,b^2}\,
\frac{D^2}{(\log n)^4\, \sqrt{n}} - \frac{c_1\,(1 + \sigma_4^2)}{n}.
$$
}

\vskip2mm
These estimates may be employed to arrive at the two-sided bounds of the form
\be
\frac{c_0}{n}\,\leq\, \E_\theta\, \omega^2(F_\theta,\Phi) 
 \,\leq\, \frac{c_1}{n}, \qquad
\frac{c_0}{(\log n)^4 \sqrt{n}} \leq 
\E_\theta\, \rho(F_\theta,\Phi) \leq \frac{c_1\log n}{\sqrt{n}}
\en
with some absolute constants $c_0 > 0$ and $c_1 > 0$.
Examples where both inequalities in (1.8) are fulfilled include
the following uniformly bounded orthonormal systems in 
$L^2(\Omega, {\mathfrak F},\P)$:

(i) The trigonometric system $X = (X_1,\dots,X_n)$ with components
\bee
X_{2k-1}(t) & = &
\sqrt{2}\,\cos(kt), \\ 
X_{2k}(t) & = &
\sqrt{2}\,\sin(kt) \quad (-\pi < t < \pi, \ k = 1,\dots,n/2, \ n \ {\rm even})
\ene
on the interval $\Omega = (-\pi,\pi)$ equipped with the normalized 
Lebesgue measure $\P$. 

(ii) The cosine trigonometric system 
$X = (X_1,\dots,X_n)$ with
$$
X_k(t) = \sqrt{2}\,\cos(kt)
$$ 
on the interval $\Omega = (0,\pi)$ equipped 
with the normalized Lebesgue measure $\P$. 

(iii) The normalized Chebyshev polynomials $X_1,\dots,X_n$
defined by
\bee
X_k(t) 
 & \hskip-2mm= &
\hskip-2mm \sqrt{2} \cos(k\arccos t)  \\
 & \hskip-2mm= &
\hskip-2mm \sqrt{2}\,\Big[t^n - 
\Big(
\begin{array}{c}
\mbox{$n$}\\
\mbox{$2$}\\ 
\end{array}
\Big) \, 
t^{n-2} (1-t^2) + 
\Big(
\begin{array}{c}
\mbox{$n$}\\
\mbox{$4$}\\ 
\end{array}
\Big) \, 
t^{n-4} (1-t^2)^2 - \dots\Big]
\ene
on $\Omega = (-1,1)$ equipped with the probability measure 
$d\P(t) = \frac{1}{\pi \sqrt{1 - t^2}}\,dt$, $|t|<1$.

(iv) The systems of functions of the form 
$$
X_k(t,s) = \Psi(kt + s), \quad k = 1,\dots,n \ \ (0 < t,s < 1)
$$
on the square $\Omega = (0,1) \times (0,1)$ equipped with the Lebesgue 
measure $\P$. In this case, (1.8) holds true for any 1-periodic Lipschitz 
function $\Psi$ on the real line such that 
$\int_0^1 \Psi(x)\,dx = 0$ and $\int_0^1 \Psi(x)^2\,dx = 1$
with constants $c_0$ and $c_1$ depending on $\Psi$ only.

(v) The Walsh system 
$$
X = \{X_\tau\}_{\tau \neq \emptyset}, \quad \tau \subset \{1,\dots,d\},
$$
of dimension $n = 2^d - 1$ on the discrete cube $\Omega = \{-1,1\}^d$ 
(the ordering of the components does not play any role). Here, $\P$ denotes
the normalized counting measure, and
$$
X_\tau(t) = \prod_{k \in \tau} t_k \quad {\rm for} \ 
t = (t_1,\dots,t_d) \in \Omega.
$$

(vi) Random vectors $X$ with associated
empirical distribution functions $F_\theta$ based on the ``observations" 
$X_k = \sqrt{n}\,\theta_k$ ($1 \leq k \leq n$).

\vskip2mm
The paper is organized as follows. We start in Section 2 with a review
of several results on the so-called typical distributions $F$ which serve as 
main approximations for $F_\theta$ (in general, they do not need to
be normal, or even nearly normal). Sections 3-7 deal
with the $L^2$-distances $\omega(F_\theta,F)$ only, while
Sections 8-12 are mostly focused on the Kolmogorov distances 
$\rho(F_\theta,F)$. In Section 13, the examples described
in items (i)-(vi) illustrate the applicability of Theorems 1.1-1.3, thus
with a standard rate of normal approximation.
In Section 14 we consider lacunary trigonometric systems and show
that the typical rate is improved to the order $1/n$.
Similar improved rates are also reviewed
in the last section in presence of certain correlation-type conditions.
Thus an outline of all sections reads as:

\vskip2mm
1. Introduction

2. Typical distributions

3. Upper bound for the $L^2$-distance at standard rate

4. General approximations for the $L^2$-distance with error 
of order at most $1/n$

5. Proof of Theorem 1.1 for the $L^2$-distance

6. General lower bounds for the $L^2$-distance. Proof of Theorem 1.2

7. Lipschitz systems

8. Berry-Esseen-type bounds

9. Quantitative forms of Sudakov's theorem for the Kolmogorov
distance

10. Proof of Theorem 1.1 for the Kolmogorov Distance

11. Relations between $L^1$, $L^2$ and Kolmogorov distances

12. Lower bounds. Proof of Theorem 1.3

13. Functional examples

14. The Walsh system; Empirical measures

15. Improved rates for lacunary systems

16. Improved rates for independent and log-concave summands

17. Improved rates under correlation-type conditions

\vskip2mm
As usual, the Euclidean space $\R^n$ is endowed with the canonical norm 
$|\,\cdot\,|$ and the inner product $\left<\cdot,\cdot\right>$. In the
sequel,  we denote by $\E_\theta$ an integral over $\S^{n-1}$ with respect 
to the measure $\s_{n-1}$. By $c$, $c_1, c_2,\dots$,  we denote positive absolute 
constants which may vary from place to place (if not stated explicitly that 
$c$ depends on some parameter). Similarly $C$ will denote a quantity bounded 
by an absolute constant. Throughout, we assume that $X$ is a given 
random vector in $\R^n$ ($n \geq 2$) and $Y$ is its independent copy.


\vskip5mm
\section{{\bf Typical Distributions}}
\setcounter{equation}{0}

\noindent
In the sequel, we denote by
$$
F(x) = \E_\theta F_\theta(x) = \E_\theta\, \P\{S_\theta \leq x\},
\quad x \in \R,
$$
the mean distribution function of the weighted sums 
$S_\theta = \left<X,\theta\right>$ with respect to the uniform measure 
$\s_{n-1}$. It is also called a typical distribution function using 
the terminology of \cite{Su}. Indeed, according to Sudakov's theorem,
if $X$ is isotropic, then most of $F_\theta$ are concentrated about $F$ 
in a weak sense (cf. \cite{A-B-P}, \cite{B1}, \cite{B-C-G3} 
for quantitative statements).

However, whether or not $F$ itself is close to the normal distribution 
function $\Phi$ is determined by the concentration properties of 
the distribution of $|X|$. Note that, due to the rotational invariance 
of $\s_{n-1}$, the typical distribution can be described as the distribution 
of the product $\theta_1\,|X|$, assuming that 
$\theta = (\theta_1,\dots,\theta_n)$ is a random vector which 
is independent of $X$ and has distribution $\s_{n-1}$. In this product, 
$\theta_1 \sqrt{n}$ is almost standard normal, so that $F$ is almost 
standard normal, if and only if $\frac{1}{\sqrt{n}}\,|X|$ is almost 1
(like in the weak law of large numbers). This assertion can be quantified 
in terms of the weighted total variation distance by virtue of 
the following upper bound derived in \cite{B-C-G2}.

\vskip5mm
{\bf Proposition 2.1.} {\sl If $\E\,|X|^2 = n$ (in particular, when $X$ is 
isotropic), then
$$
\int_{-\infty}^\infty (1 + x^2)\,|F(dx) - \Phi(dx)| \, \leq \, 
\frac{c}{n}\,\big(1 + \Var(|X|)\big).
$$
}

\vskip2mm
In particular, this gives a non-uniform bound for the normal approximation, 
namely
\be
|F(x) - \Phi(x)| \, \leq \, 
\frac{c}{n\,(1 + x^2)}\,\big(1 + \Var(|X|)\big), \quad x \in \R.
\en
In these bounds we shall rely on the following monotone functionals (of $p$)
\be
\sigma_{2p} = \sqrt{n}\,\bigg(\E\,\Big|\frac{|X|^2}{n} - 1\Big|^p\bigg)^{1/p},
\quad p \geq 1,
\en
where the particular cases $p = 1$ and $p=2$ will be most important. 
If $\E\,|X|^2 = n$, we thus deal with a more tractable quantity
$$
\sigma_4^2 = \frac{1}{n}\,\Var\big(|X|^2\big).
$$
Using an elementary inequality $\Var(\xi)\, \E \xi^2 \leq \Var(\xi^2)$
(which is true for any random random variable $\xi \geq 0$), we have
$\Var(|X|) \leq \sigma_4^2$. Another similar relation
$$
\frac{1}{4}\,\sigma_2^2 \leq \Var(|X|) \leq \sqrt{n}\,\sigma_2
$$
can be found in \cite{B-C-G3}. From (2.1), we therefore obtain the following 
bounds for the normal approximation in all $L^p$-norms
$$
\|F - \Phi\|_p = \Big(\int_{-\infty}^\infty |F(x) - \Phi(x)|^p \, dx\Big)^{1/p},
$$
including the limit case 
$$
\|F - \Phi\|_\infty = \rho(F,\Phi) = \sup_x |F(x) - \Phi(x)|.
$$

\vskip5mm
{\bf Corollary 2.2.} {\sl If $\E\,|X|^2 = n$, then, for all $p \geq 1$,
\be
\|F - \Phi\|_p \, \leq \, c\,\frac{1 + \sigma_2}{\sqrt{n}}, \qquad
\|F - \Phi\|_p \, \leq \, c\,\frac{1 + \sigma_4^2}{n}.
\en
}

\vskip2mm
Note that the  characteristic function associated to $F$ is given by
\be
f(t) \, = \, \E_\theta\, \E\,e^{it\left<X,\theta\right>} \, = \,
\E_\theta\, \E\,e^{it |X|\, \theta_1} \, = \, \E\,J_n(t|X|), \qquad t \in \R,
\en
where $J_n$ denotes the characteristic function of the first coordinate
$\theta_1$ of $\theta$ under $\s_{n-1}$. Hence, by the Plancherel theorem, 
\be
\omega^2(F,\Phi) = \frac{1}{2\pi} 
\int_{-\infty}^\infty \big(\E\,J_n(t|X|) - e^{-t^2/2}\big)^2\ \frac{dt}{t^2}.
\en
For $p=2$, the relations in (2.3) can also be derived by means of (2.5) and by 
virtue of the following Edgeworth-type approximations derived in 
\cite{B-C-G3} and \cite{B-C-G5}.

\vskip5mm
{\bf Lemma 2.3.} {\sl For all $t \in \R$,
\be
\big|J_n\big(t\sqrt{n}\big) - e^{-t^2/2}\big| \, \leq \,  
\frac{c}{n}\,\min\{1,t^2\}.
\en
Moreover,
\be
\Big|J_n\big(t\sqrt{n}) - \Big(1 - \frac{t^4}{4n}\Big)\, e^{-t^2/2}\Big|
 \, \leq \,
\frac{c}{n^2}\, \min\{1,t^4\}.
\en
}

\vskip2mm
The functions $J_n$ have a subgaussian (although oscillatory) decay 
on a long interval of the real line. In particular,  as was shown in \cite{B-C-G3}, 
\be
\big|J_n\big(t\sqrt{n}\big)\big| \, \leq \,
5\,e^{-t^2/2} + 4\,e^{-n/12}, \quad t \in \R.
\en
This bound can be used for the estimation of the characteristic function 
of the typical distribution, by involving the variance-type functionals 
$\sigma_{2p}$.

\vskip5mm
{\bf Lemma 2.4.} {\sl The characteristic function 
of the typical distribution satisfies, for all $t \in \R$,
$$
c_p\,|f(t)| \, \leq \, 
e^{-t^2/4} + \frac{1 + \sigma_{2p}^p}{n^{p/2}}
$$
with constants $c_p>0$ depending on $p \geq 1$ only.
Consequently, for all $T > 0$,
$$
\frac{c_p}{T} \int_0^T |f(t)|\,dt \, \leq \, 
\frac{1}{T} + \frac{1 + \sigma_{2p}^p}{n^{p/2}}.
$$
}

\vskip2mm
{\bf Proof.}
One may split the expectation in (2.4) to the event 
$A = \{|X|^2 \leq \lambda n\}$ and its complement 
$B = \{|X|^2 > \lambda n\}$, $0 < \lambda < 1$. By (2.8),
\bee
\E\,|J_n(t|X|)|\,1_B 
 & \leq & 
\E\,\big(5\,e^{-t^2 |X|^2/2n} + 4\,e^{-n/12}\big)\, 1_B \\
 & \leq & 
5\,e^{-\lambda t^2/2} + 4\,e^{-n/12}.
\ene
On the other hand, recalling the definition (2.2), we have
\begin{eqnarray}
\P(A) 
 & = &
\P\big\{n - |X|^2 \geq (1 - \lambda) n\big\} \nonumber \\
 & \leq & 
\frac{1}{((1 - \lambda) n)^p}\, \E\,|n - |X|^2|^p\, = \,
\frac{\sigma_{2p}^p}{(1-\lambda)^p\,n^{p/2}}.
\end{eqnarray}
Choosing $\lambda = \frac{1}{2}$, and since $|J_n(s)| \leq 1$ 
for all $s \in \R$, we get
$$
\E\,|J_n(t|X|)|\,1_A \, \leq \, (2\sigma_{2p})^p\,n^{-p/2},
$$
thus implying that
$$
|f(t)| \, \leq \, 5\, e^{-t^2/4} + 4\,e^{-n/12} + 
(2 \sigma_{2p})^p\,n^{-p/2}.
$$
This readily yields the desired pointwise and integral bounds
of the lemma.
\qed

\vskip5mm
If $|X| = \sqrt{n}$ a.s., the typical distribution $F$ is just
the distribution of $\sqrt{n}\, \theta_1$, the normalized first
coordinate of a point on the unit sphere under $\s_{n-1}$, whose 
characteristic function is $J_n(t\sqrt{n})$. In this case, the subgaussian 
character of $F$ manifests itself in corresponding deviation 
and moment inequalities such as the following.

\vskip5mm
{\bf Lemma 2.5.} {\sl For all $p > 0$,
\be
\E_\theta\, |\theta_1|^p \, \leq \, 2\,\Big(\frac{p}{n}\Big)^{p/2}.
\en
}

\vskip2mm
This inequality can be derived from the well-known
bound on the Laplace transform 
$$
\E_\theta\, e^{t\theta_1} \, \leq \, \exp\Big\{\frac{t^2}{2(n-1)}\Big\}, 
\quad t \in \R,
$$ 
which follows from the fact that the logarithmic Sobolev constant 
for the unit sphere is equal to $n-1$ (cf. \cite{L2}). Using
$x^p \leq (\frac{p}{e})^p\,e^x$, $x \geq 0$, we have 
$|x|^p \leq 2\,(\frac{p}{e})^p\,\cosh(x)$, $x \in \R$, and the above bound
implies
$$
t^p\, \E_\theta\, |\theta_1|^p \leq 
2\,\Big(\frac{p}{e}\Big)^p\, e^{\frac{t^2}{2(n-1)}}
\quad {\rm for \ all} \ t \geq 0. 
$$
The latter can be optimized over $t$, 
which leads to (2.10), even in a sharper form.

In this connection, let us emphasize that rates for the 
normal approximation for $F$ that are better than $1/n$ 
cannot be obtained under the support assumption as above.

\vskip5mm
{\bf Proposition 2.6.} {\sl For any random vector $X$ in 
$\R^n$ such that $|X|^2 = n$ a.s., we have
$$
\E_\theta\,\rho(F,\Phi)\,\geq\, \frac{c}{n}.
$$
}

\vskip2mm
{\bf Proof.} One may apply the following lower bound
\be
\rho(F,\Phi) \, \geq \, \frac{1}{3T} \ 
\Big|\int_0^T (f(t) - e^{-t^2/2})\,\Big(1 - \frac{t}{T}\Big)\, dt\Big|,
\en
which holds for any $T > 0$ (cf. \cite{B2}). Since $|X|^2 = n$ a.s., 
we have $f(t) = J_n(t\sqrt{n})$. Choosing $T=1$ and applying (2.7), 
it follows from (2.11) that $\rho(F,\Phi) \geq \frac{c}{n}$ for all 
$n \geq n_0$ where $n_0$ is determined by $c$ only.
But, a similar bound also holds for $n<n_0$ since $F$
is supported on the interval $[-\sqrt{n},\sqrt{n}]$.
\qed


\vskip5mm
\section{{\bf Upper Bound for the $L^2$-distance at Standard Rate}}
\setcounter{equation}{0}

\vskip2mm
\noindent
Like in the problem of normal approximation for the typical distribution 
function $F = \E_\theta F_\theta$, the closeness of distribution functions 
$F_\theta$ of  the weighted sums $S_\theta = \left<X,\theta\right>$ 
($\theta \in \S^{n-1}$) to $F$ in the metric $\omega$ 
can also be explored in terms of the associated characteristic functions 
(the Fourier-Stieltjes transforms)
\be
f_\theta(t) \, = \, \E\,e^{it\left<X,\theta\right>} \ = \ 
\int_{-\infty}^\infty e^{it\left<x,\theta\right>}\,dF_\theta(x), 
\quad t \in \R.
\en
Again, let us start with the identity
\be
\omega^2(F_\theta,F) = \frac{1}{2\pi} 
\int_{-\infty}^\infty \frac{|f_\theta(t)-f(t)|^2}{t^2}\,dt.
\en
Here, the mean value of the numerator represents the variance
$\E_\theta\, |f_\theta(t)|^2 - |f(t)|^2$ with respect to $\s_{n-1}$. 
Moreover, using an independent copy $Y$ of $X$, we have
\be
\E_\theta\, |f_\theta(t)|^2 \, = \, \E_\theta\,
\E\, e^{it\left<X-Y,\theta\right>} \, = \, \E J_n(t|X-Y|).
\en
Hence, the Plancherel formula (3.2) together with (2.4) yields
\be
\E_\theta\, \omega^2(F_\theta,F) = \frac{1}{2\pi} 
\int_{-\infty}^\infty 
\Big(\E J_n(t|X-Y|) - \big(\E J_n(t|X|)\big)^2\Big)\,\frac{dt}{t^2}.
\en

In this section our aim is to show that the above expression is 
of order at most $O(1/n)$ provided that the mean $a = \E X$,
$m_2 = m_2(X)$ and $\sigma_4^2 = \sigma_4^2(X)$ are of order 1. 
The next statement contains the upper bound (1.2) as
a partial case.

\vskip5mm
{\bf Proposition 3.1.} {\sl Given a random vector $X$ in $\R^n$ 
with $\E X = a$ and $\E\,|X|^2 = n$, we have
\be
\E_\theta\, \omega^2(F_\theta,F) \, \leq \, \frac{cA}{n}
\en
with $A = 1+ |a|^2 + m_2^2  + \sigma_4^2$.
A similar inequality continues to hold with the
normal distribution function $\Phi$ in place of $F$.
}

\vskip5mm
If $X$ is isotropic, then $m_2 = 1$, while $|a| \leq 1$ (by Bessel's
inequality). Hence, both characteristics $m_2$ and $a$ may be 
removed from the parameter $A$ in this case. However,
in the general case, it may happen that $m_2$ and $\sigma_4$
are bounded, while $|a|$ is large. The example in Remark 3.2
shows that this parameter can not be removed.

\vskip5mm
{\bf Proof.} Note that, for any $\eta > 0$,
\be
\int_{-\infty}^\infty \frac{\min\{1,t^2 \eta^2\}}{t^2}\,dt \, = \, 4\eta,
\en
Hence,  in the formula (3.4), the expectation $\E J_n(t|X-Y|)$ 
can be replaced using the normal approximation (2.6) 
at the expense of an error not exceeding
$$
\frac{c}{n}\, \E \int_{-\infty}^\infty 
\min\Big\{1,\frac{t^2 |X-Y|^2}{n}\Big\}\,\frac{dt}{t^2} \ = \
\frac{4c}{n}\, \E \, \frac{|X-Y|}{\sqrt{n}} \ \leq \ \frac{8c}{n},
$$
where we used that $\E \,|X| \leq \sqrt{n}$.
Similarly, by (2.6) and (3.6),
\bee
\int_{-\infty}^\infty \Big|\,
\big(\E J_n(t|X|)\big)^2 - \big(\E\,e^{-t^2|X|^2/2n}\big)^2 \Big|\,
\frac{dt}{t^2} 
 & \leq &
2\, \E \int_{-\infty}^\infty 
\big|J_n(t|X|) - e^{-t^2|X|^2/2n}\big|\,\frac{dt}{t^2} \\
 & \leq &
\frac{2c}{n}\, \E \int_{-\infty}^\infty 
\min\Big\{1,\frac{t^2 |X|^2}{n}\Big\}\,\frac{dt}{t^2} \\
 & = &
\frac{8c}{n}\, \E \, \frac{|X|}{\sqrt{n}} \ \leq \ \frac{8c}{n}.
\ene
Hence, using these bounds in (3.4), we arrive at the general 
approximation
\be
\E_\theta\, \omega^2(F_\theta,F) \, = \, \frac{1}{2\pi} 
\int_{-\infty}^\infty \Big(\E\,e^{-t^2|X-Y|^2/2n} - 
\big(\E\, e^{-t^2 |X|^2/2n}\big)^2\Big)\,\frac{dt}{t^2} + 
\frac{C}{n},
\en
where we recall that $C$ denotes a quantity bounded by 
an absolute constant.

Introduce the random variable
$$
\rho^2 = \frac{|X-Y|^2}{2n} \quad (\rho \geq 0).
$$
By Jensen's inequality, $\E\,e^{-t^2|X|^2/2n} \geq e^{-t^2/2}$,
so that, by (3.7),
$$
\E_\theta\, \omega^2(F_\theta,F) \, \leq \,
\frac{1}{2\pi}\,\E \int_{-\infty}^\infty 
\frac{e^{-\rho^2 t^2} - e^{-t^2}}{t^2}\,dt + \frac{c}{n}.
$$
The above integral is easily evaluated (by differentiating
with respect to the variable ``$\rho^2$"), and we arrive at the bound
\be
\E_\theta\, \omega^2(F_\theta,F) \leq \frac{1}{\sqrt{\pi}}\,
(1 - \E \rho) + \frac{c}{n}.
\en

To further simplify, one may apply an elementary inequality
$1 - x\, \leq \, \frac{1}{2}\,(1-x^2) + (1-x^2)^2$ ($x \geq 0$),
which gives
$$
\E_\theta\, \omega^2(F_\theta,F) \leq 
\frac{1}{2\sqrt{\pi}}\,\E\,(1 - \rho^2) +
\frac{1}{\sqrt{\pi}}\, \E\,(1 - \rho^2)^2 + \frac{c}{n}.
$$
Since
$$
1 - \rho^2 = 
\frac{n - |X|^2}{2n} + \frac{n - |Y|^2}{2n} + 
\frac{\left<X,Y\right>}{n},
$$
we have
$$
1 - \E \rho^2 = \frac{1}{n}\,\E\left<X,Y\right> = 
\frac{1}{n}\,|\E X|^2 = \frac{1}{n}\,|a|^2.
$$
In addition,
$$
(1 - \rho^2)^2 \, \leq \,
2\,\bigg(\frac{n - |X|^2}{2n} + \frac{n - |Y|^2}{2n}\bigg)^2 + 
2\,\frac{\left<X,Y\right>^2}{n^2},
$$
which implies
$$
\E\,(1 - \rho^2)^2 \, \leq \,
\frac{\Var(|X|^2)}{n^2} + 2\,\frac{\E \left<X,Y\right>^2}{n^2} 
 \, = \, \frac{\sigma_4^2 + 2m_2^2}{n}.
$$
Using this estimate in (3.8), the inequality (3.5) follows
immediately.

For the second assertion, it remains to apply Corollary 2.2.
\qed

\vskip5mm
{\bf Remark 3.2.} Let us illustrate the inequality (3.5) in the example 
where the random vector $X$ has a normal distribution with a large 
mean value. Given a standard normal random vector 
$Z = (Z_1,\dots,Z_{n-1})$ in $\R^{n-1}$ (which we identify with 
the space of all points in $\R^n$ with zero last coordinate), define
$$
X = \alpha Z + \lambda e_n \quad {\rm with} \ \ 
1 \leq \lambda \leq n^{1/4}, \ \ \alpha^2 (n-1) + \lambda^2 = n,
$$
where $e_n = (0,\dots,0,1)$ is the last unit vector in the
canonical basis of $\R^n$. Since $Z$ is orthogonal to $e_n$,
so that $|X|^2 = \alpha^2\,|Z|^2 + \lambda^2$, we have
$\E\,|X|^2 = n$, and
$$
\sigma_4^2 \, = \, \frac{\alpha^4}{n}\,\Var(|Z|^2) \, = \,
\frac{2\alpha^4\,(n-1)}{n} \, = \, 2\,\frac{(n-\lambda^2)^2}{n(n-1)}
\, < \, 2.
$$

Let $Z'$ be an independent copy of $Z$. Then
$Y = \alpha Z' + \lambda e_n$ is an independent copy of $X$, so that
$$
m_2^2 \, = \, \frac{1}{n}\,\E \left<X,Y\right>^2 \, = \,
\frac{1}{n}\,(\alpha^4\,(n-1) + \lambda^4) \, < \, 2.
$$
Thus, both $m_2$ and $\sigma_4$ are bounded, while the mean 
$a = \E X = \lambda e_n$ has the Euclidean length 
$|a| = \lambda \geq 1$. Hence, the inequality (3.5) being stated
for the normal distribution function in place of $F$ simplifies to
$$
\E_\theta\, \omega^2(F_\theta,\Phi) \, \leq \, \frac{c\lambda^2}{n}.
$$

Let us show that this bound may be reversed up to an absolute
factor (which would imply that $|a|^2$ may not be removed from $A$). 
For any unit vector $\theta = (\theta_1,\dots,\theta_n)$, the linear form 
$$
S_\theta = \left<X,\theta\right> = 
\alpha \theta_1 Z_1 + \dots + \alpha \theta_{n-1} Z_{n-1} + \lambda \theta_n
$$ 
has a normal distribution on the line with mean 
$\E S_\theta = \lambda \theta_n$ and variance 
$\Var(S_\theta) = \alpha^2 (1 - \theta_n^2)$.
Consider the normal distribution function 
$\Phi_{\mu,\sigma^2}(x) = \Phi(\frac{x - \mu}{\sigma})$ with parameters 
$0 \leq \mu \leq 1$ and $\frac{1}{2} \leq \sigma^2 \leq 1$
($\sigma>0$). If $x \leq \frac{\mu}{1+\sigma}$, then 
$\frac{x - \mu}{\sigma} \leq x$, and on the interval with these 
endpoints the standard normal density $\varphi(y)$ attains minimum 
at the left endpoint. Hence
$$
|\Phi_{\mu,\sigma^2}(x) - \Phi(x)| \, = \,
\int_{\frac{x - \mu}{\sigma}}^x \varphi(y)\,dy \, \geq \,
\Big(x - \frac{x - \mu}{\sigma}\Big)\,\varphi\Big(\frac{x - \mu}{\sigma}\Big),
$$
so that
\bee
\omega^2(\Phi_{\mu,\sigma^2},\Phi) 
 & \geq & 
\int_{-\infty}^{\frac{\mu}{1+\sigma}} 
\Big(x - \frac{x - \mu}{\sigma}\Big)^2\,
\varphi\Big(\frac{x - \mu}{\sigma}\Big)^2\,dx \\
 & = &
\frac{\sigma}{2\pi}  \int_{-\infty}^{-\frac{\mu}{1+\sigma}} 
(\mu - (1-\sigma) y)^2\,e^{-y^2/2}\,dy\, \geq \,
\frac{\sigma \mu^2}{2\pi}\,\int_{-\infty}^{-\frac{\mu}{1+\sigma}} 
e^{-y^2/2}\,dy \, \geq \, c\mu^2.
\ene
In our case, since $\lambda \leq n^{1/4}$ and
$$
\alpha^2 \, = \, \frac{n - \lambda^2}{n-1} \, \geq \, 
\frac{n - \sqrt{n}}{n-1} \, \geq \, 1 - \frac{1}{\sqrt{n}},
$$ 
we have $|\E S_\theta| \leq 1$ and
$\Var(S_\theta) \geq \frac{1}{2}$ on the set 
$\Omega_n = \{\theta \in \S^{n-1}: |\theta_n| < \frac{\log n}{\sqrt{n}}\}$
with $n$ large enough. It follows that 
$$
\E_\theta\, \omega^2(F_\theta,\Phi)  \, \geq \, 
c\lambda^2\,\E\,\theta_n^2\,1_{\{\theta \in \Omega_n\}}
 \, \geq \, \frac{c' \lambda^2}{n}.
$$


\vskip5mm
\section{{\bf General Approximations for the $L^2$-distance with Error 
of Order at most $1/n$}}
\setcounter{equation}{0}

\vskip2mm
\noindent
We now turn to general representations for the average $L^2$-distance
between $F_\theta$ and the typical distribution function $F$ with
error of order at most $1/n$.

\vskip5mm
{\bf Proposition 4.1.} {\sl Suppose that $\E\,|X| \leq b\sqrt{n}$
for some $b \geq 0$. Then
\be
\E_\theta\, \omega^2(F_\theta,F) = 
\frac{1}{\sqrt{2\pi}}\,\E R + \frac{Cb}{n^2},
\en
where
\be
R \, = \, \frac{(|X|^2 + |Y|^2)^{1/2}}{\sqrt{n}}\,
\Big(1 + \frac{1}{4n}\,
\frac{|X|^4 + |Y|^4}{(|X|^2 + |Y|^2)^2}\Big) - 
\frac{|X-Y|}{\sqrt{n}}\, \Big(1 + \frac{1}{4n}\Big).
\en
}

\vskip5mm
We use the convention that $R=0$ if $X=Y=0$.
Note that $|R| \leq 3\,\frac{|X| + |Y|}{\sqrt{n}}$, so $\E R \leq 3b$.

Let us give a simpler expression by involving the functional
$\sigma_4^2 = \frac{1}{n}\,\Var(|X|^2)$ and 
assuming that $\E\,|X|^2 = n$. Since
$$
\frac{|X|^4 + |Y|^4}{(|X|^2 + |Y|^2)^2} - \frac{1}{2} 
 \, = \,
\frac{(|X|^2 - |Y|^2)^2}{2\,(|X|^2 + |Y|^2)^2},
$$
we may write
\be
R \, = \, \frac{1}{8n^{3/2}}\,
\frac{(|X|^2 - |Y|^2)^2}{(|X|^2 + |Y|^2)^{3/2}} +
\frac{(|X|^2 + |Y|^2)^{1/2}}{\sqrt{n}}\,\Big(1 + \frac{1}{8n}\Big) - 
\frac{|X-Y|}{\sqrt{n}}\, \Big(1 + \frac{1}{4n}\Big).
\en
As we will see, the first term here is actually of order at most $\sigma_4^2/n^2$.
As a result, we arrive at the relation (1.5).

\vskip5mm
{\bf Proposition 4.2.} {\sl If $\E\,|X|^2 = n$, then
\be
\E_\theta\, \omega^2(F_\theta,F) = 
\frac{1}{\sqrt{2\pi}}\,\E R + C\,\frac{1 + \sigma_4^2}{n^2},
\en
where
\be
R \, = \, \frac{(|X|^2 + |Y|^2)^{1/2}}{\sqrt{n}}\,
\Big(1 + \frac{1}{8n}\Big) - 
\frac{|X-Y|}{\sqrt{n}}\, \Big(1 + \frac{1}{4n}\Big).
\en
}

\vskip2mm
{\bf Proof of Proposition 4.1.} Let us return to the Plancherel formula
(3.4). To simplify the integrand therein, we apply the inequality (2.7) 
in Lemma 2.3, by replacing $t^4$ with $t^2$ in the remainder term. 
Using the equality (3.6), the expectation $\E J_n(t|X-Y|)$ in the formula 
(3.4) can be therefore replaced according to (2.7) at the expense of 
an error not exceeding
$$
\frac{c}{n^2}\, \E \int_{-\infty}^\infty 
\min\Big\{1,\frac{t^2 |X-Y|^2}{n}\Big\}\,\frac{dt}{t^2} \ = \
\frac{4c}{n^2}\, \E \, \frac{|X-Y|}{\sqrt{n}} \ \leq \ \frac{8cb}{n^2}.
$$

As for the main term $(1 - \frac{t^4}{4n})\, e^{-t^2/2}$ 
in (2.7), it is bounded by an absolute constant, which implies that
\bee
J_n\big(t\sqrt{n}) J_n\big(s\sqrt{n}) 
 & = &
\Big(1 - \frac{t^4}{4n}\Big)\Big(1 - \frac{s^4}{4n}\Big)\, 
e^{-(t^2 + s^2)/2} + O\big(n^{-2} \min\{1,t^2 + s^2\}\big) \\ 
 & = &
\Big(1 - \frac{t^4 + s^4}{4n}\Big)\, 
e^{-(t^2 + s^2)/2} + O\big(n^{-2} \min\{1,t^2 + s^2\}\big).
\ene
Hence
\bee
|\E J_n(t|X|)|^2 
 \, = \,
\E\, J_n(t|X|) \, J_n(t|Y|)
 & = &
\E\, \Big(1 - \frac{t^4\, (|X|^4 + |Y|^4)}{4n^3}\Big)\, 
e^{-\frac{t^2\, (|X|^2 + |Y|^2)}{2n}} \\
 & & + \  
O\Big(n^{-2} \min\Big\{1,\frac{t^2\, (|X|^2 + |Y|^2)}{n}\Big\}\Big).
\ene
As before, after integration in (3.4) the latter remainder term will produce
a quantity not exceeding a multiple of $b/n^2$. As a preliminary step,
we therefore obtain the representation
\be
\E_\theta\, \omega^2(F_\theta,F) = \frac{1}{2\pi}\,I + \frac{Cb}{n^2}
\en
with
$$
I \, = \, \E \int_{-\infty}^\infty  
\bigg[\Big(1 - \frac{t^4 |X-Y|^4}{4n^3}\Big)\,
e^{-\frac{t^2|X-Y|^2}{2n}} - 
\Big(1 - \frac{t^4\, (|X|^4 + |Y|^4)}{4n^3}\Big)\, 
e^{-\frac{t^2\, (|X|^2 + |Y|^2)}{2n}}\bigg]\,\frac{dt}{t^2}.
$$

To evaluate the integrals of this type, consider the functions
$$
\psi_r(\alpha) = \frac{1}{\sqrt{2\pi}} \int_{-\infty}^\infty 
\Big((1 - rt^4)\,e^{-\alpha t^2/2} - e^{-t^2/2}\Big)\,\frac{dt}{t^2}
 \qquad (\alpha > 0, \ r \in \R).
$$
Clearly,
$$
\psi_r(1) = -\frac{1}{\sqrt{2\pi}}  
\int_{-\infty}^\infty rt^2\,e^{-t^2/2}\,dt = -r
$$
and
\bee
\psi_r'(\alpha) 
 & = &
-\frac{1}{2\sqrt{2\pi}} \int_{-\infty}^\infty 
(1 - rt^4)\,e^{-\alpha t^2/2}\,dt \\
 & = &
-\frac{1}{2\sqrt{\alpha}}\, \frac{1}{\sqrt{2\pi}} 
\int_{-\infty}^\infty 
\Big(1 - \frac{r}{\alpha^2}\,s^4\Big)\,e^{-s^2/2}\,ds
 \, = \, 
-\frac{1}{2\sqrt{\alpha}}\,\Big(1 - \frac{3r}{\alpha^2}\Big).
\ene
Hence
\bee
\psi_r(\alpha) - \psi_r(1)
 & = &
\int_1^\alpha \Big(-\frac{1}{2}\,z^{-1/2} + 
\frac{3r}{2}\,z^{-5/2}\Big)\,dz \\
 & = &
(1+r) - (\alpha^{1/2} + r\alpha^{-3/2}),
\ene
and we get
\be
\psi_r(\alpha) = 1 - (\alpha^{1/2} + r\alpha^{-3/2}).
\en
Here, when $\alpha$ and $r$ both approach zero subject to the relation
$r = O(\alpha^2)$, we get in the limit $\psi_0(0) = 1$. From this,
\bee
\frac{1}{\sqrt{2\pi}}\,I 
 & = & 
\E\,(\psi_{r_1}(\alpha_1) - \psi_{r_2}(\alpha_2)) \\
 & = &
\E\,(\alpha_2^{1/2} + r_2\alpha_2^{-3/2}) - 
\E\,(\alpha_1^{1/2} + r_1\alpha_1^{-3/2}),
\ene
which we need with
\bee
\alpha_1 = \frac{|X-Y|^2}{n}, \ \
 & & 
r_1 = \frac{|X-Y|^4}{4n^3}, \\
\alpha_2 = \frac{|X|^2 + |Y|^2}{n}, \ \ 
 & &
r_2 = \frac{|X|^4 + |Y|^4}{4n^3}.
\ene
It follows that
\bee
\alpha_2^{1/2} + r_2\alpha_2^{-3/2} 
 & = &
\Big(\frac{|X|^2 + |Y|^2}{n}\Big)^{1/2} \,
\Big(1 + \frac{1}{4n}\,\frac{|X|^4 + |Y|^4}{(|X|^2 + |Y|^2)^2}\Big), \\
\alpha_1^{1/2} + r_1\alpha_1^{-3/2}
 & = &
\bigg(\frac{|X-Y|^2}{n}\bigg)^{1/2} \Big(1 + \frac{1}{4n}\Big)
\ene
with the assumption that both expressions are equal to zero 
in the case $X=Y=0$. As a result, (4.6) yields the desired
representation (4.1) with quantity $R$ described in (4.2).
\qed

\vskip5mm
In order to modify (4.1)-(4.2) to the form (4.4)-(4.5), first let us verify
the following general relation.

\vskip5mm
{\bf Lemma 4.3.} {\sl Let $\xi$ be a non-negative random variable 
with finite second moment (not identically zero), and let $\eta$ be 
its independent copy. Then
$$
\E\,\frac{(\xi - \eta)^2}{(\xi + \eta)^{3/2}}\,1_{\{\xi + \eta > 0\}} 
 \, \leq \, 12\,\frac{\Var(\xi)}{(\E\,\xi)^{3/2}}.
$$
}

\vskip2mm
Applying the lemma with $\xi = |X|^2$, $\eta = |Y|^2$ and 
assuming that $\E\,|X|^2 = n$, we get that
$$
\E\,
\frac{(|X|^2 - |Y|^2)^2}{(|X|^2 + |Y|^2)^{3/2}} \, \leq \,
12\,\frac{\Var(|X|^2)}{(\E\,|X|^2)^{3/2}} \, = \,  
12\,\frac{\Var(|X|^2)}{n^{3/2}} \, = \,
12\,\frac{\sigma_4^2}{n^{1/2}}.
$$
In view of (4.3), this proves Proposition 4.2.

\vskip5mm
{\bf Proof of Lemma 4.3.} By homogeneity, we may assume that 
$\E \xi = 1$. In particular, $\E\,|\xi - \eta| \leq 2$. We have
\bee
\E\,\frac{(\xi - \eta)^2}{(\xi + \eta)^{3/2}}\,
1_{\{\xi + \eta > 1/2\}} 
 & \leq &
2^{3/2} \,\E\,(\xi - \eta)^2\,1_{\{\xi + \eta > 1/2\}} \\
 & \leq &
2^{3/2} \,\E\,(\xi - \eta)^2 \, = \, 4 \sqrt{2}\,\Var(\xi).
\ene
Also note that, by Chebyshev's inequality,
$$
\P\,\{\xi \leq 1/2\} \, = \, \P\,\{1 - \xi \geq 1/2\} \, \leq \,
4\,\Var(\xi)^2,
$$
so
$$
\P\,\{\xi + \eta \leq 1/2\} \, \leq \,
\P\,\{\xi \leq 1/2\}\, \P\,\{\eta \leq 1/2\} \, \leq \, 16\,\Var(\xi)^2.
$$
Hence, since $\frac{|\xi - \eta|}{\xi + \eta} \leq 1$
for $\xi + \eta > 0$, we have, by Cauchy's inequality,
\bee
\E\,\frac{(\xi - \eta)^2}{(\xi + \eta)^{3/2}}\,
1_{\{0 < \xi + \eta \leq 1/2\}}
 & \leq &
\E\, \sqrt{|\xi - \eta|}\,1_{\{0 < \xi + \eta \leq 1/2\}} \\
 & \leq &
\sqrt{\E\,|\xi - \eta|}\, \sqrt{\P\,\{\xi + \eta \leq 1/2\}}
 \, \leq \, 4\sqrt{2}\,\Var(\xi).
\ene
It remains to combine both inequalities, which yield
$$
\E\,\frac{(\xi - \eta)^2}{(\xi + \eta)^{3/2}}\,1_{\{\xi + \eta > 0\}} 
\, \leq \, 8\sqrt{2}\,\Var(\xi) \, \leq \, 12\,\Var(\xi).
$$
\qed


\vskip5mm
\section{{\bf Proof of Theorem 1.1 for the $L^2$-distance}}
\setcounter{equation}{0}

\vskip2mm
\noindent
The expression (4.5) may be further simplified in the particular case where
the distribution of $X$ is supported on the sphere $\sqrt{n}\ \S^{n-1}$.
Introduce the random variable
$$
\xi = \frac{\left<X,Y\right>}{n},
$$ 
where $Y$ is an independent copy of $X$. Since $|X-Y|^2 = 2n\,(1 - \xi)$,
Proposition 4.2 yields:

\vskip5mm
{\bf Corollary 5.1.} {\sl If $|X|^2 = n$ a.s., then
\be
\sqrt{\pi}\ \E_\theta\, \omega^2(F_\theta,F)  \, = \,
\Big(1 + \frac{1}{4n}\Big)\,\E\,\Big(1 - (1-\xi)^{1/2}\Big) - 
\frac{1}{8n} + O\Big(\frac{1}{n^2}\Big).
\en
}

\vskip2mm
Note that $|\xi| \leq 1$. Therefore, the relation (5.1) suggests
to develop an expansion in powers of $\ep$ for the function 
$w(\ep) = 1 - \sqrt{1 - \ep}$ near zero, which will be needed up to 
the term $\ep^4$.

\vskip5mm
{\bf Lemma 5.2.} {\sl For all $|\ep| \leq 1$,
$$
1 - \sqrt{1 - \ep} \ \leq \ \frac{1}{2}\,\ep + \frac{1}{8}\,\ep^2 +
\frac{1}{16}\,\ep^3 + 3\ep^4.
$$
In addition,
$$
1 - \sqrt{1 - \ep} \ \geq \ \frac{1}{2}\,\ep + \frac{1}{8}\,\ep^2 +
\frac{1}{16}\,\ep^3 + 0.01\,\ep^4.
$$
}

\vskip2mm
{\bf Proof.} By Taylor's formula for the function $w(\ep)$ 
around zero on the half-axis $\ep < 1$, 
$$
1 - \sqrt{1 - \ep} \ = \ \frac{1}{2}\,\ep + \frac{1}{8}\,\ep^2 +
\frac{1}{16}\,\ep^3 + \frac{5}{128}\,\ep^4 + 
\frac{w^{(5)}(\ep_1)}{120}\,\ep^5
$$
for some $\ep_1$ between zero and $\ep$. Since
$w^{(5)}(\ep) = \frac{105}{32}\,(1-\ep)^{-9/2} \geq 0$, we have 
an upper bound
$$
1 - \sqrt{1 - \ep} \, \leq \, \frac{1}{2}\,\ep + \frac{1}{8}\,\ep^2 +
\frac{1}{16}\,\ep^3 + \frac{5}{128}\,\ep^4, \qquad \ep \leq 0.
$$
Also, $w^{(5)}(\ep) \leq \frac{105}{32}\,3^{9/2} < 461$ for 
$0 \leq \ep \leq \frac{2}{3}$, so, in this interval
$$
\frac{5}{128}\,\ep^4 + \frac{w^{(5)}(\ep_1)}{120}\,\ep^5 \leq
3\ep^4.
$$
Thus, in both cases,
$$
1 - \sqrt{1 - \ep} \, \leq \, \frac{1}{2}\,\ep + \frac{1}{8}\,\ep^2 +
\frac{1}{16}\,\ep^3 + 3\ep^4, \qquad
\ep \leq \frac{2}{3}.
$$
To treat the remaining values $\frac{2}{3} \leq \ep \leq 1$, it 
is sufficient to select a positive constant $b$ such that the polynomial
$$
Q(\ep) = \frac{1}{2}\,\ep + \frac{1}{8}\,\ep^2 +
\frac{1}{16}\,\ep^3 + b\ep^4
$$
is greater than or equal to $1$ for $\ep \geq \frac{2}{3}$. On this
half-axis, $Q(\ep) \geq \frac{11}{27} + b\,\frac{16}{81} \geq 1$ 
for $b \geq 3$. Thus, the upper bound of the lemma is proved.

Now, from Taylor's formula we also get that
$$
1 - \sqrt{1 - \ep} \, \geq \, \frac{1}{2}\,\ep + \frac{1}{8}\,\ep^2 +
\frac{1}{16}\,\ep^3 + \frac{5}{128}\,\ep^4, \qquad \ep \geq 0.
$$
In addition, if $-1 \leq \ep \leq 0$, then 
$w^{(5)}(\ep) \leq \frac{105}{32}$, so
\bee
1 - \sqrt{1 - \ep} 
 & = & 
\frac{1}{2}\,\ep + \frac{1}{8}\,\ep^2 +
\frac{1}{16}\,\ep^3 + \frac{5}{128}\,\ep^4\, 
\Big(1 + \frac{w^{(5)}(\ep_1)}{120}\,\ep\Big) \\
 & \geq & 
\frac{1}{2}\,\ep + \frac{1}{8}\,\ep^2 +
\frac{1}{16}\,\ep^3 + \ep^4\, 
\Big(\frac{5}{128} - \frac{\frac{105}{32}}{120}\Big) \\
 & \geq & 
\frac{1}{2}\,\ep + \frac{1}{8}\,\ep^2 +
\frac{1}{16}\,\ep^3 + 0.01\, \ep^4.
\ene
\qed

\vskip2mm
{\bf Proof of Theorem 1.1} (First part).
Using Lemma 5.2 with $\ep = \xi$ and 
applying Corollary 5.1, we get an asymptotic representation
$$
\sqrt{\pi}\ \E_\theta\, \omega^2(F_\theta,F)  \, = \,
\Big(1 + \frac{1}{4n}\Big)\,
\Big(\frac{1}{8}\,\E \xi^2 + \frac{1}{16}\,\E \xi^3 + c\,\E \xi^4\Big) 
- \frac{1}{8n} + O\Big(\frac{1}{n^2}\Big)
$$
for some quantity $c$ such that $0.01 \leq c \leq 3$. 
If additionally $X$ is isotropic, then $\E \left<X,Y\right>^2 = n$, i.e.
$\E \xi^2 = \frac{1}{n}$, and the representation is simplified to
$$
\sqrt{\pi}\ \E_\theta\, \omega^2(F_\theta,F)  \, = \,
\Big(1 + \frac{1}{4n}\Big)\,
\Big(\frac{1}{16}\,\E \xi^3 + c\,\E \xi^4\Big) 
+ O\Big(\frac{1}{n^2}\Big),
$$
thus removing the term of order $1/n$. Moreover, since
$\E \xi^4 \leq \E\, |\xi|^3 \leq \E \xi^2 = \frac{1}{n}$, the fraction
$\frac{1}{4n}$ may be removed from the brackets at the expense 
of the remainder term. Thus
$$
\sqrt{\pi}\ \E_\theta\, \omega^2(F_\theta,F)  \, = \,
\frac{1}{16}\,\E \xi^3 + c\,\E \xi^4 + O\Big(\frac{1}{n^2}\Big),
$$
which is exactly the expansion (1.3).
\qed

\vskip5mm
{\bf Remark 5.3.} In the isotropic case with $|X|^2 = n$ a.s.,
but without the mean zero assumption, the above expansion takes the form
\be
\sqrt{\pi}\ \E_\theta\, \omega^2(F_\theta,F)  \, = \,
\frac{1}{2}\,\E \xi +
\frac{1}{16}\,\E \xi^3 + c\,\E \xi^4 + O\Big(\frac{1}{n^2}\Big).
\en
Since the last two expectations are non-negative, this implies in
particular that
\be
\E_\theta\, \omega^2(F_\theta,F)  \, \geq \,
\frac{1}{2\sqrt{\pi}}\,\E \xi + O\Big(\frac{1}{n^2}\Big).
\en


\vskip5mm
\section{{\bf General Lower Bounds for the $L^2$-distance. Proof of Theorem 1.2}}
\setcounter{equation}{0}

\vskip2mm
\noindent
Proposition 4.1 may be used to establish the following general lower bound 
which will be the first step in the proof of Theorem 1.2. Recall that $Y$
denotes an independent copy of a random vector $X$ in $\R^n$.

\vskip5mm
{\bf Proposition 6.1.} {\sl If $\E\,|X| \leq b\sqrt{n}$, then
\be
\E_\theta\, \omega^2(F_\theta,F) \, \geq \, 
c_1\,\E\, \rho\,\xi^4 - c_2 \frac{b}{n^2},
\en
where
$$
\rho = \Big(\frac{|X|^2 + |Y|^2}{2n}\Big)^{1/2}, \qquad
\xi = \frac{2 \left<X,Y\right>}{|X|^2 + |Y|^2}.
$$
}

\vskip2mm
The argument employs two elementary lemmas. 

\vskip5mm
{\bf Lemma 6.2.} {\sl If $\E\,|X|^2$ is finite, then
\be
\E \left<X,Y\right>^2 \, \geq \, \frac{1}{n}\,\big(\E\,|X|^2\big)^2.
\en
}

\vskip2mm
By the invariance of (6.2) under 
linear orthogonal transformations, we may assume that 
$\E X_i X_j = \lambda_i \delta_{ij}$ where $\lambda_i$'s appear as 
eigenvalues of the covariance operator of $X$. Since
$$
\E\, |X|^2 = \sum_{i=1}^n \lambda_i, \qquad
\E \left<X,Y\right>^2  = \sum_{i=1}^n \lambda_i^2,
$$
the inequality (6.2) follows by applying Cauchy's inequality.

\vskip5mm
{\bf Lemma 6.3.} {\sl If $\E\,|X|^p$ is finite for an integer
$p \geq 1$, then, for any real number $0 \leq \alpha \leq p$,
$$
\E\, \frac{\left<X,Y\right>^p}{(|X|^2 + |Y|^2)^\alpha} \geq 0,
$$
where the ratio is defined to be zero in case $X=Y=0$.
In addition, for $\alpha \in [0,2]$,
$$
\E\, \frac{\left<X,Y\right>^2}{(|X|^2 + |Y|^2)^\alpha} \, \geq \, \frac{1}{n}\,
\E\, \frac{|X|^2\, |Y|^2}{(|X|^2 + |Y|^2)^\alpha}.
$$
}

\vskip2mm
{\bf Proof.} First, let us note that
$$
\E\,\frac{|\left<X,Y\right>|^p}{(|X|^2 + |Y|^2)^\alpha} \, \leq \, 
\E\,\frac{(|X|\,|Y|)^p}{(|X|\,|Y|)^\alpha} \, = \, (\E\,|X|^{p-\alpha})^2,
$$
so, the expectation on the left is finite. Without loss of generality, 
we may assume that $0 < \alpha \leq p$ and $r = |X|^2 + |Y|^2 > 0$ 
with probability 1. We use the identity
$$
\int_0^\infty e^{-r t^{1/\alpha}} dt \, = \, c_\alpha\, r^{-\alpha} \quad 
{\rm where} \ \ c_\alpha = \int_0^\infty e^{-s^{1/\alpha}}\,ds,
$$
which gives
$$
c_\alpha\, \E \left<X,Y\right>^p r^{-\alpha} \, = \,
\int_0^\infty \E \left<X,Y\right>^p\,e^{-r t^{1/\alpha}}\,dt.
$$
Writing $X = (X_1,\dots,X_n)$ and $Y = (Y_1,\dots,Y_n)$, we have
\bee
\E\left<X,Y\right>^p\, e^{-r t^{1/\alpha}} 
 & = &
\E\left<X,Y\right>^p\, e^{-t^{1/\alpha} (|X|^2 + |Y|^2)} \\
 & = &
\sum_{i_1,\dots,i_p = 1}^n 
\Big(\E\,X_{i_1} \dots X_{i_p} \ e^{-t^{1/\alpha}\, |X|^2}\Big)^2,
\ene
which shows that the left expectation is always non-negative. 
Integrating over $t>0$, this proves the first assertion.

For the second assertion, write
$$
c_\alpha\,\E \left<X,Y\right>^2 r^{-\alpha} \, = \,
\int_0^\infty \E \left<X,Y\right>^2\, e^{-t^{1/\alpha} (|X|^2 + |Y|^2)}\,dt
 \, = \, \int_0^\infty \E \left<X_t,Y_t\right>^2\,dt,
$$
where
$$
X_t = e^{-t^{1/\alpha} |X|^2/2}\,X, \quad 
Y_t = e^{-t^{1/\alpha} |Y|^2/2}\,Y.
$$
Since $Y_t$ represents an independent copy of $X_t$, one may apply 
Lemma 6.2 which gives
$$
\E \left<X_t,Y_t\right>^2 \, \geq \, \frac{1}{n}\,\E\,|X_t|^2\,|Y_t|^2.
$$
Hence,
\bee
\int_0^\infty \E \left<X_t,Y_t\right>^2\,dt 
 & \geq & 
\frac{1}{n}\, \int_0^\infty \E\,|X_t|^2\,|Y_t|^2\,dt \\
 & = &
\frac{1}{n}\, \int_0^\infty \E\,|X|^2\,|Y|^2\, e^{-t^{1/\alpha} 
(|X|^2 + |Y|^2)}\,dt
 \ = \
\frac{c_\alpha}{n}\,\E\, |X|^2\,|Y|^2\, r^{-\alpha}.
\ene
\qed

\vskip2mm
{\bf Proof of Proposition 6.1.}
Let us return to the representation (4.3) in Proposition 4.1 and
write
$$
\E_\theta\, \omega^2(F_\theta,F) = 
\frac{1}{\sqrt{2\pi}}\,\E\, (R_0 + R_1) + \frac{Cb}{n^2},
$$
where
$$
R_0 \, = \, \frac{1}{8n^{3/2}}\,
\frac{(|X|^2 - |Y|^2)^2}{(|X|^2 + |Y|^2)^{3/2}}
$$
and
\bee
R_1 
 & = &
\frac{(|X|^2 + |Y|^2)^{1/2}}{\sqrt{n}}\,\Big(1 + \frac{1}{8n}\Big) - 
\frac{|X-Y|}{\sqrt{n}}\, \Big(1 + \frac{1}{4n}\Big) \\
 & = &
\frac{(|X|^2 + |Y|^2)^{1/2}}{\sqrt{n}}\,\bigg[\Big(1 + \frac{1}{4n}\Big) 
\big(1 - \sqrt{1 - \xi}\,\big) - \frac{1}{8n}\bigg]
\ene
with the assumption that $R_0 = 0$ when $X=Y=0$. Since $|\xi| \leq 1$, one
may apply Lemma 5.2 which gives
$$
R_1 \ \geq \ \frac{(|X|^2 + |Y|^2)^{1/2}}{\sqrt{n}}\,
\bigg[\Big(1 + \frac{1}{4n}\Big) 
\Big(\frac{1}{2}\,\xi + \frac{1}{8}\,\xi^2 +
\frac{1}{16}\,\xi^3 + 0.01\, \xi^4\Big) - \frac{1}{8n}\bigg].
$$
The expectation of the terms on the right-hand side containing 
$\xi$ and $\xi^3$ is non-negative according to Lemma 6.3 with 
$\alpha = \frac{1}{2}$, $p = 1$, and with $\alpha = \frac{5}{2}$, 
$p = 3$, respectively. Hence, removing the unnecessary factor 
$1 + \frac{1}{4n}$, we get
\begin{eqnarray}
\E_\theta\, \omega^2(F_\theta,F) 
 & \geq &
\frac{1}{\sqrt{2\pi}}\,\E R_0 + \frac{1}{\sqrt{2\pi}}\,
\E\, \frac{(|X|^2 + |Y|^2)^{1/2}}{8\sqrt{n}}\, 
\Big(\xi^2 - \frac{1}{n}\Big) \nonumber \\
 & & + \ 
c_1\, \E\, \frac{(|X|^2 + |Y|^2)^{1/2}}{\sqrt{n}}\,\xi^4 -
c_2 \frac{b}{n^2}.
\end{eqnarray}

Now, by the second inequality of Lemma 6.3 applied with
$\alpha = 3/2$, $p = 2$, we have
\bee
\E\, (|X|^2 + |Y|^2)^{1/2}\,\xi^2 
 & = &
4\ \E\, \frac{\left<X,Y\right>^2}{(|X|^2 + |Y|^2)^{3/2}} \\
 & \geq &
\frac{4}{n}\,\E\, \frac{|X|^2\,|Y|^2}{(|X|^2 + |Y|^2)^{3/2}}.
\ene
This gives
\bee
\E\, \frac{(|X|^2 + |Y|^2)^{1/2}}{8\sqrt{n}}\,
\Big(\xi^2 - \frac{1}{n}\Big) 
 & \geq &
\frac{1}{8n^{3/2}}\,\E\, 
\Big[\frac{4\,|X|^2\,|Y|^2}{(|X|^2 + |Y|^2)^{3/2}} - 
(|X|^2 + |Y|^2)^{1/2}\Big] \\
 & = &
- \frac{1}{8n^{3/2}}\,
\E\, \frac{(|X|^2 - |Y|^2)^2}{(|X|^2 + |Y|^2)^{3/2}} \, = \,
- \E R_0.
\ene
Thus, the  summand $\E R_0$ in (6.3) neutralizes the second expectation, 
and we are left with the term containing $\xi^4$. 
\qed

\vskip5mm
{\bf Proof of Theorem 1.2.} We apply Proposition 6.1.
By the assumption, $\E \rho^2 = 1$ and
$\Var(\rho^2) = \frac{1}{2n}\,\sigma_4^2$, where 
$\sigma_4^2 = \frac{1}{n}\,\Var(|X|^2)$. Using
$$
2 \left<X,Y\right> = |X|^2 + |Y|^2 - |X-Y|^2, \qquad
\xi = 1 - \frac{|X-Y|^2}{|X|^2 + |Y|^2},
$$
we have
\bee
\xi^4 
 & \geq &
(1-\alpha)^4\,
1_{\{|X-Y|^2 \, \leq \, \alpha\,(|X|^2 + |Y|^2)\}} \\
 & \geq &
(1-\alpha)^4\,
1_{\{|X-Y|^2 \, \leq \, \alpha \lambda n, \ 
|X|^2 + |Y|^2 \, \geq \, \lambda n\}}, 
\quad 0 < \alpha, \lambda < 1.
\ene
On the set $|X|^2 + |Y|^2 \geq \lambda n$, we necessarily
have $\rho^2 \geq \frac{\lambda}{2}$, so
\bee
\E\, \rho\,\xi^4 
 & \geq &
\frac{(1-\alpha)^4}{\sqrt{2}}\, \sqrt{\lambda} \ 
\P\Big\{|X-Y|^2 \leq \alpha \lambda n, \ 
|X|^2 + |Y|^2 \geq \lambda n\Big\} \\ 
 & \geq &
\frac{(1-\alpha)^4}{\sqrt{2}}\, \sqrt{\lambda} \ 
\Big(\P\{|X-Y|^2 \leq \alpha \lambda n\} - 
\P\{|X|^2 + |Y|^2 \leq \lambda n\}\Big).
\ene
But, by Chebyshev's inequality 
$$
\P\big\{|X|^2 \leq \lambda n\big\} \, = \,
\P\big\{n - |X|^2 \geq (1-\lambda)\, n\big\} \, \leq \,
\frac{\Var(|X|^2)}{(1-\lambda)^2\, n^2} \, = \,  
\frac{\sigma_4^2}{(1-\lambda)^2\, n},
$$
implying
$$
\P\big\{|X|^2 + |Y|^2 \leq \lambda n\big\} \, \leq \,
\Big(\P\big\{|X|^2 \leq \lambda n\big\}\Big)^2  \, \leq \,
\frac{1}{(1-\lambda)^4} \ \frac{\sigma_4^4}{n^2}.
$$ 
Hence
$$
\E\, \rho\,\xi^4 \, \geq \,
\frac{(1-\alpha)^4}{\sqrt{2}}\, \sqrt{\lambda} \ 
\bigg(\P\{|X-Y|^2 \leq \alpha \lambda n\} - 
\frac{1}{(1-\lambda)^4} \ 
\frac{\sigma_4^4}{n^2}\bigg).
$$
Choosing, for example, $\alpha = \lambda = \frac{1}{2}$, 
we get
$$
\E\, \rho\,\xi^4 \, \geq \,
\frac{1}{32}\,\P\Big\{|X-Y|^2 \leq \frac{1}{4}\, n\Big\} - 
\frac{\sigma_4^4}{2n^2}.
$$
It remains to apply (6.1) with $b=1$ and replace $F$ with $\Phi$
on the basis of (2.3).
\qed


\vskip5mm
\section{{\bf Lipschitz Systems}}
\setcounter{equation}{0}

\vskip2mm
\noindent
While upper bounds of order $n^{-1/2}$ for the $L^2$-distance 
$\omega(F_\theta,F)$ on average are provided in (1.2) and in 
the more general inequality (3.5) of Proposition 3.1, in this section 
we focus on the conditions that provide similar 
lower bounds, as a consequence of Theorem 1.2.

Let $L$ be a fixed measurable function on the underlying probability 
space $(\Omega,{\mathfrak F},\P)$. We will say that the system 
$X_1,\dots,X_n$ of random variables on $(\Omega,{\mathfrak F},\P)$, 
or the random vector $X = (X_1,\dots,X_n)$ in $\R^n$ satisfies 
a Lipschitz condition with a parameter function $L$, if
\be
\max_{1 \leq k \leq n} |X_k(t) - X_k(s)| \leq n\,|L(t) - L(s)|, \qquad 
t,s \in \Omega.
\en
When $\Omega$ is an interval of the real line (finite or not), and 
$L(t) = Lt$, $L > 0$, this condition means that every function 
$X_k$ in the system has a Lipschitz semi-norm at most $Ln$.

As before, we use the variance functional 
$\sigma_4^2 = \frac{1}{n}\,\Var(|X|^2)$.

\vskip5mm
{\bf Proposition 7.1.} {\sl Suppose that $\E\,|X|^2 = n$. 
If the random vector $X$ satisfies the Lipschitz condition with 
a parameter function $L$, then
\be
\E_\theta\, \omega^2(F_\theta,F) \,\geq\, 
\frac{c_L}{n} - \frac{c_0\,(1 + \sigma_4^4)}{n^2}
\en
with some absolute constant $c_0>0$ and with a constant $c_L$ 
depending on the distribution of $L$ only. Moreover, if $L$ has 
finite second moment, then with some absolute constant $c_1>0$
\be
\E_\theta\, \omega^2(F_\theta,F) \,\geq\, 
\frac{c_1}{n\sqrt{\Var(L)}} - \frac{c_0\, (1 + \sigma_4^4)}{n^2}.
\en
}

\vskip2mm
Note that, if $X_1,\dots,X_n$ form an orthonormal system 
in $L^2(\Omega,{\mathfrak F},\P)$, i.e., the random vector $X$
is isotropic, and if $L$ has finite second moment 
$\|L\|_2^2 = \E L^2$, then this moment has to be bounded from 
below by a multiple of $1/n^2$. Indeed, the projection of 
the function $\eta(t) = 1$ in  $L^2(\Omega,{\mathfrak F},\P)$ 
to the linear hull $H$ of $X_1,\dots,X_n$ has the form
${\rm Proj}_H(\eta) = \sum_{k=1}^n \left<\eta,X_k\right> X_k$,
and we have Bessel's inequality
$$
1 = \|\eta\|_2^2 \geq \|{\rm Proj}_H(\eta)\|_2^2 =
\sum_{k=1}^n \left<\eta,X_k\right>^2 = \sum_{k=1}^n\, (\E X_k)^2
$$
(where we used the canonical innde product $\left<\cdot,\cdot\right>$
in $L^2(\Omega,{\mathfrak F},\P)$). By the Lipschitz assumption,
$|X_k(t) - X_k(s)|^2 \leq n^2\,|L(t) - L(s)|^2$. Integrating this 
inequality over the product measure $\P(dt) \otimes \P(ds)$, we obtain 
a lower bound
$$
n^2\, \Var(L) \geq \Var(X_k) = 1 - (\E X_k)^2.
$$
One may now perform summation over $k = 1,\dots, n$, which 
together with Bessel's inequality leads to
$$
\Var(L) \geq \frac{n-1}{n^3} \geq \frac{1}{2n^2} \quad (n \geq 2).
$$

The Lipschitz condition (7.1) guarantees the validity of 
the following property, which can be combined with Theorem 1.2
to obtain (7.2)-(7.3).

\vskip5mm
{\bf Lemma 7.2.} {\sl Suppose that the random vector 
$X = (X_1,\dots,X_n)$ satisfies the Lipschitz condition with 
the parameter function $L$. If $Y$ is 
an independent copy of $X$, then
$$
\P\big\{|X - Y|^2 \leq \lambda n\big\} \geq 
\frac{c \sqrt{\lambda}}{n}, \qquad 0 \leq \lambda \leq 1,
$$
where the constant $c>0$ depends on the distribution of 
$L$ only. Moreover, if $L$ has finite second moment, then
$$
\P\big\{|X - Y|^2 \leq \lambda n\big\} \geq 
\frac{\sqrt{\lambda}}{6n\sqrt{\Var(L)}}, \qquad
0 \leq \lambda \leq n^2\, \Var(L).
$$
}

\vskip5mm
In turn, this lemma is based on the following general 
observation.

\vskip5mm
{\bf Lemma 7.3.} {\sl If $\eta$ is an independent copy of 
a random variable $\xi$, then for any $\ep_0 > 0$, 
$$
\P\{|\xi - \eta| \leq \ep\} \geq 
c\ep, \qquad 0 \leq \ep \leq \ep_0,
$$
with some constant $c>0$ independent of $\ep$. Moreover, if 
the standard deviation $\sigma = \sqrt{\Var(\xi)}$ is finite, then
$$
\P\{|\xi - \eta| \leq \ep\} \geq \frac{1}{6\sigma}\,\ep, \qquad 
0 \leq \ep \leq \sigma.
$$
}

\vskip2mm
{\bf Proof.} The difference $\xi - \eta$ has a non-negative characteristic 
function $h(t) = |\psi(t)|^2$, where $\psi$ is the characteristic function 
of $\xi$. Denoting by $H$ the distribution function of
$\xi - \eta$, we start with a general identity
\be
\int_{-\infty}^\infty \hat p(x)\,dH(x) = 
\int_{-\infty}^\infty p(t)h(t)\,dt,
\en
which is valid for any integrable function $p(t)$ on the real 
line with Fourier transform 
$\hat p(x) = \int_{-\infty}^\infty e^{itx}\,p(t)\,dt$, $x \in \R$.
Given $\ep>0$, here we take a standard pair
$$
p(t) = \frac{1}{2\pi}\,
\Big(\frac{\sin \frac{\ep t}{2}}{\frac{\ep t}{2}}\Big)^2, 
\qquad \hat p(x) = \frac{1}{\ep}\,\Big(1 - \frac{|x|}{\ep}\Big)^+,
$$
where we use the notation $a^+ = \max\{a,0\}$. In this case, 
$$
\int_{-\infty}^\infty \hat p(x)\,dH(x) \, \leq \,
\frac{1}{\ep} \int_{[-\ep,\ep]} dH(x) \, = \, 
\frac{1}{\ep}\,\P\{|\xi - \eta| \leq \ep\}.
$$

On the other hand, since the function $\frac{\sin u}{u}$ is 
decreasing in $0 < u < \frac{\pi}{2}$, we have
$$
\int_{-\infty}^\infty p(t)h(t)\,dt \, \geq \, 
\frac{1}{2\pi}\,\big(2\sin(1/2)\big)^2
\int_{-1/\ep}^{1/\ep} h(t)\,dt\, \geq \, 
\frac{1}{7} \int_{-1/\ep}^{1/\ep} h(t)\,dt.
$$
Hence, whenever $0 < \ep \leq \ep_0$, by (7.4),
$$
\P\{|\xi - \eta| \leq \ep\} \, \geq \, 
\frac{\ep}{7} \int_{-1/\ep}^{1/\ep} h(t)\,dt\, \geq \, 
\frac{\ep}{7} \int_{-1/\ep_0}^{1/\ep_0} h(t)\,dt.
$$
Since $h(t)$ is bounded away from zero near the origin, 
the first assertion follows. 

One may quantify this statement in terms of the variance 
$\sigma^2 = \Var(\xi)$ by using Taylor's expansion for 
$h(t)$ about zero. Indeed, it gives $1 - h(t) \leq \sigma^2 t^2$, 
and thus for $\ep \leq \ep_0 = \sigma$, 
$$ 
\int_{-1/\ep}^{1/\ep} h(t)\,dt \, \geq \,
\int_{-1/\sigma}^{1/\sigma} (1 - \sigma^2 t^2)\,dt 
 \, = \,\frac{4}{3\sigma}.
$$
Since $\frac{\ep}{7} \cdot \frac{4}{3\sigma} \geq 
\frac{1}{6\sigma}\,\ep$, the lemma is proved.
\qed

\vskip5mm
{\bf Proof of Lemma 7.2.} Let us equip the product space 
$\Omega^2 = \Omega \times \Omega$ with the product measure 
$\P^2 = \P \otimes \P$ and redefine $X$ on this new probability 
space as $X(t,s) = X(t)$, $(t,s) \in \Omega^2$. Then one 
can introduce an independent copy of $X$ in the form 
$Y(t,s) = X(s)$. By the Lipschitz condition,
$$
|X(t,s) - Y(t,s)|^2 \, = \, 
\sum_{k=1}^n |X_k(t) - X_k(s)|^2 \, \leq \, n^3 \,|L(t) - L(s)|^2.
$$
Hence, if $\eta$ is an independent copy of the random variable 
$\xi = L$, then
$$
\P\big\{|X - Y|^2 \leq \lambda n\big\} \, \geq \, 
\P\big\{n^3\, |\xi - \eta|^2 \leq \lambda n\big\} \, = \,
\P\bigg\{|\xi - \eta| \leq \frac{\sqrt{\lambda}}{n}\bigg\}.
$$
But, by Lemma 7.3 with $\ep_0 = 1$, the latter probability is 
at least $c\,\frac{\sqrt{\lambda}}{n}$, where the constant $c$ 
depends on $L$ only (via its distribution). An application of the 
second inequality of Lemma 7.3 yields the second assertion.
\qed

\vskip5mm
To include more examples, let us now give a bit more general 
form of Lemma 7.2, assuming that 
$(\Omega,\P) = (\Omega_1 \times \Omega_2,\P_1 \otimes \P_2)$
is a product probability space.

\vskip5mm
{\bf Lemma 7.4.} {\sl Let 
$X = (X_1,\dots,X_n): \Omega \rightarrow \R^n$ be a random 
vector such that, for some measurable functions $L_1$ and $L_2$ 
defined on $\Omega_1$ and $\Omega_2$ respectively, 
\be
\max_{1 \leq k \leq n} |X_k(t_1,t_2) - X_k(s_1,s_2)| \, \leq \, 
n\,|L_1(t_1) - L_1(s_1)| + |L_2(t_2) - L_2(s_2)| 
\en
for all $(t_1,t_2), (s_1,s_2) \in \Omega$. If $Y$ is an independent 
copy of $X$, then
\be
\P\big\{|X - Y|^2 \leq \lambda n\big\} \geq \frac{c\lambda}{n}, 
 \qquad 0 \leq \lambda \leq 1,
\en
where the constant $c>0$ depends on the distributions of 
$L_1$ and $L_2$ only.
}

\vskip5mm
{\bf Proof.} Again, let us equip the product space 
$\Omega^2 = \Omega \times \Omega$ with the product measure 
$\P^2 = \P \otimes \P$ and put $X(t,s) = X(t)$, $Y(t,s) = X(s)$
for $t = (t_1,t_2) \in \Omega$ and $s = (s_1,s_2) \in \Omega$,
so that $Y$ is an independent copy of $X$.
By the Lipschitz condition (7.5), for any $k \leq n$,
$$
|X_k(t) - X_k(s)|^2  \, \leq \, 
2n^2\,|L_1(t_1) - L_1(s_1)| + 2\,|L_2(t_2) - L_2(s_2)|^2,
$$
so
\bee
|X(t) - Y(s)|^2 
 & = &
\sum_{k=1}^n |X_k(t) - X_k(s)|^2 \\
 & \leq &
2n^3\,|L_1(t_1) - L_1(s_1)|^2 + 2n\,|L_2(t_2) - L_2(s_2)|^2.
\ene
Putting $L_1(t_1,t_2) = L_1(t_1)$ and $L_2(t_1,t_2) = L_2(t_2)$, 
one may treat $L_1$ and $L_2$ as independent random variables. 
If $L_1'$ is an independent copy of $L_1$ and $L_2'$ is 
an independent copy of $L_2$, we obtain that
\bee
\P\big\{|X - Y|^2 \leq \lambda n\big\} 
 & \geq &
\P\Big\{n^2\, |L_1 - L_1'|^2 + |L_2 - L_2'|^2 \leq 
\frac{\lambda}{2}\Big\} \\
 & \geq &
\P\Big\{n^2\, |L_1 - L_1'|^2 \leq \frac{\lambda}{4}\Big\}\,
\P\Big\{|L_2 - L_2'|^2 \leq \frac{\lambda}{4}\Big\} \\
 & = &
\P\Big\{|L_1 - L_1'| \leq \frac{1}{2n}\sqrt{\lambda}\,\Big\}\,
\P\Big\{|L_2 - L_2'| \leq \frac{1}{2}\sqrt{\lambda}\,\Big\}.
\ene
It remains to apply Lemma 7.3.
\qed

\vskip2mm
Let us now combine the inequality (1.8) of Theorem 1.2 with
the inequality (7.6) applied with $\lambda = \frac{1}{4}$. 
Then we obtain the following generalization of Proposition 7.1.

\vskip5mm
{\bf Proposition 7.5.} {\sl Under the Lipschitz condition 
$(7.5)$, we have
$$
\E_\theta\, \omega^2(F_\theta,F) \,\geq\, 
\frac{c}{n} - \frac{c_0\,(1 + \sigma_4^4)}{n^2},
$$
where $c_0>0$ is an absolute constant, while
$c>0$ depends on the distributions of 
$L_1$ and $L_2$. A similar estimate also holds when 
$F$ is replaced with the normal distribution function~$\Phi$.
}

\vskip5mm
The last assertion follows from the inequality
(2.3), cf. Corollary 2.2.


\vskip10mm
\section{{\bf Berry-Esseen-type Bounds}}
\setcounter{equation}{0}

\vskip2mm
\noindent
We now turn to the study of the Kolmogorov distance 
$$
\rho(F_\theta,F) = \sup_x\, |F_\theta(x) - F(x)|, \quad
\theta \in \S^{n-1},
$$
between the distribution functions $F_\theta$ of the weighted sums
$S_\theta = \left<X,\theta\right>$ and the typical distribution function
$F = \E_\theta F_\theta$. We are mostly interested in bounding 
the second moment $\E_\theta\, \rho^2(F_\theta,F)$.
As in the case of the $L^2$-distance, our basic tool will be a Fourier 
analytic approach relying upon a general Berry-Esseen-type bound
\be
c\,\rho(U,V) \, \leq \,
\int_0^T \frac{|\hat U(t) - \hat V(t)|}{t}\,dt + \frac{1}{T}  
\int_0^T |\hat V(t)|\,dt, \qquad T>0,
\en
where $U$ and $V$ may be arbitrary distribution functions on the line  
with characteristic functions $\hat U$ and $\hat V$ respectively 
(cf. e.g. \cite{B2}, \cite{P1}, \cite{P2}). 

As before, we denote by $f_\theta$ and $f$ the characteristic functions 
associated to $F_\theta$ and $F$. Recall that $\sigma_{2p}$-functionals 
were defined in (2.2).

\vskip5mm
{\bf Lemma 8.1.} {\sl If $T \geq T_0 \geq 1$, then for all $p \geq 1$,
\begin{eqnarray}
c_p\,\E_\theta\,\rho^2(F_\theta,F) 
 & \leq &
\int_0^1 \frac{\E_\theta\,|f_\theta(t) - f(t)|^2}{t^2}\,dt +
\log T \, \int_0^{T_0} \frac{\E_\theta\,|f_\theta(t) - f(t)|^2}{t}\,dt 
\nonumber \\
 & & \hskip-2mm + \ 
\log T \, \int_{T_0}^T \frac{\E_\theta\,|f_\theta(t)|^2}{t}\,dt +
\frac{1}{T^2} + \frac{1 + \sigma_{2p}^{2p}}{n^p},
\end{eqnarray}
where the constants $c_p>0$ depend on $p$ only.
}

\vskip5mm
{\bf Proof.} By (8.1), for any $\theta \in \S^{n-1}$,
$$
c\,\rho(F_\theta,F) \, \leq \, \int_0^T 
\frac{|f_\theta(t) - f(t)|}{t}\,dt + \frac{1}{T} \int_0^T |f(t)|\,dt,
$$
and squaring it, we get
$$
c\, \rho^2(F_\theta,F) \, \leq \,
\Big(\int_0^T \frac{|f_\theta(t) - f(t)|}{t}\,dt\Big)^2 + 
\frac{1}{T^2}\, \Big(\int_0^T |f(t)|\,dt\Big)^2.
$$
Let us split integration in the first integral into the intervals $[0,1]$ 
and $[1,T]$. By Cauchy's inequality,
$$
\Big(\int_0^1 \frac{|f_\theta(t) - f(t)|}{t}\,dt\Big)^2 \, \leq \,  
\int_0^1 \frac{|f_\theta(t) - f(t)|^2}{t^2}\,dt,
$$
while
$$
\Big(\int_1^T \frac{|f_\theta(t) - f(t)|}{t}\,dt\Big)^2 \, \leq \,
\log T \, \int_1^T \frac{|f_\theta(t) - f(t)|^2}{t}\,dt.
$$
Hence
\bee
c\,\rho^2(F_\theta,F) 
 & \leq &
\int_0^1 \frac{|f_\theta(t) - f(t)|^2}{t^2}\,dt \\
 & & + \ 
\log T \, \int_1^T \frac{|f_\theta(t) - f(t)|^2}{t}\,dt +
\frac{1}{T^2}\, \Big(\int_0^T |f(t)|\,dt\Big)^2.
\ene

Without an essential loss one may extend integration in 
the second integral to the larger interval $[0,T]$. Moreover, 
taking the expectation over $\theta$, we then get
\bee
c\,\E_\theta\,\rho^2(F_\theta,F) 
 & \leq &
\int_0^1 \frac{\E_\theta\,|f_\theta(t) - f(t)|^2}{t^2}\,dt \\
 & & + \ 
\log T \, \int_0^T \frac{\E_\theta\,|f_\theta(t) - f(t)|^2}{t}\,dt +
\frac{1}{T^2}\, \Big(\int_0^T |f(t)|\,dt\Big)^2.
\ene

Again, one may split integration in the second last integral to the two 
intervals $[0,T_0]$ and $[T_0,T]$, so that to consider separately 
sufficiently large values of $t$ for which $|f_\theta(t)|$ is small enough 
(with high probability). More precisely, since
$f(t) = \E_\theta\,f_\theta(t)$ and
$$
|f_\theta(t) - f(t)|^2 \leq 2\,|f_\theta(t)|^2 + 2\,|f(t)|^2,
$$
we have $|f(t)|^2 \leq \E_\theta\,|f_\theta(t)|^2$ and therefore
$$
\E_\theta\,|f_\theta(t) - f(t)|^2 \leq 4\,\E_\theta\,|f_\theta(t)|^2.
$$
It remains to apply Lemma~2.4.
\qed

\vskip5mm
In order to control the last integral in (8.2), one may apply the upper bound 
(2.8) on $J_n$ in the representation (3.3) to get that, for all $t \in \R$,
$$
\E_\theta\, |f_\theta(t)|^2\, \leq \, 
5\,\E\,e^{-t^2 |X-Y|^2/2n} + 4\,e^{-n/12},
$$
where $Y$ is an independent copy of the random vector $X$.
Splitting the last expectation to the event
$A = \{|X-Y|^2 \leq \frac{1}{4}\, n\}$ and its complement leads~to 
\be
\E_\theta\, |f_\theta(t)|^2 \, \leq \, 
5\,e^{-t^2/8} + 4\,e^{-n/12} + 5\,\P(A).
\en
The latter probability may further be estimated by using the moment
functionals such as $m_p$. 

To recall the argument (cf. also \cite{B-C-G3}, Proposition 2.5), 
first note that, by (2.9) with $\lambda = \frac{3}{4}$,
$$
\P\Big\{|X|^2 + |Y|^2 \leq \frac{3}{4}\, n\Big\} \, \leq \,
\P\Big\{|X|^2 \leq \frac{3}{4}\, n\Big\}\,  
\P\Big\{|Y|^2 \leq \frac{3}{4}\, n\Big\} \, \leq \,
\frac{(4\sigma_{2p})^{2p}}{n^p}.
$$
On the other hand, by Markov's inequality, assuming that $p \geq 1$ is integer,
we have
$$
\P\Big\{|\left<X,Y\right>| \geq \frac{1}{4}\, n\Big\} \, \leq \,
\frac{4^{2p}\, \E \left<X,Y\right>^{2p}}{n^{2p}} \, = \,
\frac{4^{2p}\,m_{2p}^{2p}}{n^p}.
$$
Since $|X-Y|^2 = |X|^2 + |Y|^2 - 2\left<X,Y\right>$, we have
$$
\Big\{|X-Y|^2 \leq \frac{1}{4}\,\Big\} \subset
\Big\{|X|+|Y|^2 \leq \frac{1}{4}\,n\Big\} \cup 
\Big\{\left<X,Y\right> > \frac{1}{4}\,n\Big\},
$$
and it follows that
$$
\P(A) \, \leq \, \P\Big\{|X|^2 + |Y|^2 \leq \frac{3}{4}\, n\Big\} +
\P\Big\{\left<X,Y\right> > \frac{1}{4}\, n\Big\} \, \leq \, 
\frac{4^{2p}}{n^p}\,(m_{2p}^{2p} + \sigma_{2p}^{2p}).
$$

Returning to (8.3) and noting that necessarily 
$m_{2p} \geq m_2 \geq 1$ under the assumption that 
$\E\,|X|^2 = n$, we thus obtain that
$$
c_p\,\E_\theta\, |f_\theta(t)|^2 \, \leq \,
\frac{m_{2p}^{2p} + \sigma_{2p}^{2p}}{n^p} + e^{-t^2/8}.
$$
Using this bound, the inequality (8.2) is simplified:

\vskip5mm
{\bf Lemma 8.2.} {\sl If the random vector $X$ in $\R^n$ satisfies 
$\E\,|X|^2 = n$, then for all $T \geq T_0 \geq 1$ and any integer $p \geq 1$,
\begin{eqnarray}
c_p\,\E_\theta\,\rho^2(F_\theta,F) 
 & \leq &
\int_0^1 \frac{\E_\theta\,|f_\theta(t) - f(t)|^2}{t^2}\,dt +
\log T \, \int_0^{T_0} \frac{\E_\theta\,|f_\theta(t) - f(t)|^2}{t}\,dt 
\nonumber \\
 & & \hskip-2mm + \ 
\frac{m_{2p}^{2p} + \sigma_{2p}^{2p}}{n^p}\, (1 + \log T)^2 + 
\frac{1}{T^2} + e^{-T_0^2/8}\,\log T
\end{eqnarray}
with constants $c_p$ depending on $p$ only.
}


\vskip10mm
\section{{\bf Quantitative Forms of Sudakov's Theorem for the Kolmogorov
Distance}}
\setcounter{equation}{0}

\vskip2mm
\noindent
Let us specialize Lemma 8.2 to the value $p=1$, assuming that the random vector 
$X$ is isotropic in $\R^n$ (so that $m_2 = 1$). If $\sigma_2$ is bounded, then 
choosing 
$$
T = 4n, \quad T_0 = 4\sqrt{\log n},
$$
the last three terms in (8.4) produce a quantity of order at most 
$(\log n)^2/n$. In order to bound the integrals in (8.4), one may apply 
the classical Poincar\'e inequality on the unit sphere $\S^{n-1}$
\be
\E_\theta |u(\theta)|^2 \, \leq \, 
\frac{1}{n-1}\,\E_\theta\,|\nabla u(\theta)|^2
\en
to the mean zero functions $u_t(\theta) = f_\theta(t) - f(t)$. They are 
well defined and smooth on $\R^n$ for any fixed value $t \in \R$ and have
gradients (by differentiating in (3.1)) given by
$$
\left<\nabla u_t(\theta),w\right> \, = \, 
it\,\E\,\left<X,w\right> e^{it\left<X,\theta\right>}, \quad w \in \C^n,
$$
where we use the canonical inner product in the product complex space.
By the isotropy assumption, 
$$
|\left<\nabla u_t(\theta),w\right>| \leq
|t|\,\E\,|\left<X,w\right>| \leq |t|\,|w|
$$ 
for all $w$. Hence 
$|\nabla u_t(\theta)|^2 \leq t^2$ for any $\theta \in \R^n$, so that by (9.1),
\be
\E_\theta\,|f_\theta(t) - f(t)|^2 \leq \frac{t^2}{n-1}.
\en
Applying this inequality in (8.4) together with the first bound
in (2.3) in order to replace $F$ with $\Phi$, we obtain:

\vskip5mm
{\bf Proposition 9.1.} {\sl Given an isotropic random vector
$X$ in $\R^n$,
$$
\E_\theta\, \rho^2(F_\theta,\Phi) \, \leq \, 
c\,(1 + \sigma_2^2)\,\frac{(\log n)^2}{n}.
$$
}

\vskip2mm
Since $\sigma_2 \leq \sigma_4$, we thus have
\be
\big(\E_\theta\, \rho^2(F_\theta,\Phi)\big)^{1/2} \, \leq \, c\,
(1 + \sigma_4)\,\frac{\log n}{\sqrt{n}}
\en
which sharpens (1.1). The latter bound will be an essential step 
in the proof of Theorem 1.3, while (1.1) is not strong enough.

Let us now consider another scenario in Lemma 8.2, where the 
distribution of $X$ is supported on the sphere $\sqrt{n}\ \S^{n-1}$. 
In this case,
\bee
\E_\theta\,|f_\theta(t) - f(t)|^2 
 & = &
\E_\theta\,|f_\theta(t)|^2 - |f(t)|^2 \\
 & = &
\E J_n(t|X-Y|) - J_n(t\sqrt{n})^2
\ene
according to (3.3), while $\sigma_4 = 0$. Hence, in (8.4) with $p=2$ 
we arrive at the following preliminary bound which is needed for 
the proof of Theorem 1.1 in its second part. Here we use again that
$m_4 \geq m_2 \geq 1$.

\vskip5mm
{\bf Corollary 9.2.} {\sl Suppose that $|X| = \sqrt{n}$ a.s., and
$Y$ is an independent copy of $X$. Then
\be
c\,\E_\theta\,\rho^2(F_\theta,F) \, \leq \,
\int_0^1 \frac{\Delta_n(t)}{t^2}\,dt + 
\log n \, \int_0^{4\sqrt{\log n}} \frac{\Delta_n(t)}{t}\,dt + 
\frac{(\log n)^2}{n^2}\,m_4^4,
\en
where
\be
\Delta_n(t) = \E J_n(t|X-Y|) - J_n(t\sqrt{n})^2.
\en
}


\vskip5mm
\section{{\bf Proof of Theorem 1.1 for the Kolmogorov Distance}}
\setcounter{equation}{0}

\vskip2mm
\noindent
To study the integrals in (9.4), assume additionally that the random vector 
$X$ in $\R^n$ is isotropic with mean zero and put
$$
\xi = \frac{\left<X,Y\right>}{n}, 
$$
where $Y$ is an independent copy of $X$. Note that 
$\frac{1}{n^2}\,m_4^4 = \E \xi^4$ which is present in the last term 
on the right-hand side of (9.4).

Focusing on the first integral, we need to develop an asymptotic 
bound on $\Delta_n(t)$ for $t \in [0,1]$. 
Since $|X-Y|^2 = 2n (1 - \xi)$, (9.5) becomes
$$
\Delta_n(t) = 
\E J_n\big(t\sqrt{2n (1-\xi)}\,\big) - \big(J_n(t\sqrt{n})\big)^2.
$$
We use the asymptotic formula (2.7),
\be
J_n\big(t\sqrt{n}) = \Big(1 - \frac{t^4}{4n}\Big)\,e^{-t^2/2} +
\ep_n(t), \quad t \in \R,
\en
where $\ep_n(t)$ denotes a quantity of the form 
$O\big(n^{-2} \min(1,t^4)\big)$ with a universal constant in $O$.
It implies a similar representation
\be
\big(J_n\big(t\sqrt{n})\big)^2 = \Big(1 - \frac{t^4}{2n}\Big)\,e^{-t^2} +
\ep_n(t).
\en
Since $|\xi| \leq 1$ a.s., we also have
$$
J_n\big(t\sqrt{2n (1-\xi)}\,\big) =
\Big(1 - \frac{t^4}{n}\, (1-\xi)^2\Big)\,e^{-t^2 (1-\xi)} + 
\ep_n(t).
$$
Hence, subtracting from $e^{t^2 \xi}$ the linear term $1 + t^2 \xi$
and adding, one may write
\bee
\Delta_n(t) 
 & = &
e^{-t^2}\,\E\, \Big(\Big(1 - \frac{t^4}{n}\, (1-\xi)^2\Big)\,e^{t^2 \xi} -
\Big(1 - \frac{t^4}{2n}\Big)\Big) + \ep_n(t) \\
 & = &
e^{-t^2}\,\E\,(U + V) + \ep_n(t)
\ene
with
\bee
U 
 & = &
\frac{t^4}{n}\,\Big(\frac{1}{2} - (1 - \xi)^2\Big) + 
\Big(1 - \frac{t^4}{n}\,(1-\xi)^2\Big) \cdot t^2 \xi, \\
V
 & = &
\Big(1 - \frac{t^4}{n}\,(1-\xi)^2\Big) (e^{t^2 \xi} - 1 - t^2 \xi).
\ene
Using $\E \xi = 0$, $\E \xi^2 = \frac{1}{n}$ and hence 
$\E\, |\xi|^3 \leq \E \xi^2 \leq \frac{1}{n}$, we find that
in the interval $0 \leq t \leq 1$,
$$
\E\, U = - \frac{t^4}{2n} - \frac{t^4}{n^2} + \frac{2t^6}{n^2} - 
\frac{t^6}{n}\,\E \xi^3 = - \frac{t^4}{2n} + \ep_n(t).
$$

Next write
$$
V = W - \frac{t^4}{n}\, (1-\xi)^2\, W, \qquad 
W = e^{t^2 \xi} - 1 - t^2 \xi.
$$
Using $|e^x - 1 - x| \leq 2x^2$ for $|x| \leq 1$, we have
$|W| \leq 2t^4 \xi^2$. Hence, the expected value of the second 
term in the representation for $V$ does not exceed $8t^8/n^2$. 
Moreover, by Taylor's expansion,
$$
W = \frac{1}{2}\,t^4\xi^2 + \frac{1}{6}\,t^6\xi^3 + R t^8 \xi^4, 
\qquad R = \sum_{k=4}^\infty \frac{t^{2k-8}}{k!}\,\xi^{k-4},
$$
implying that 
$$
\E\,W = \frac{t^4}{2n} + \frac{t^6}{6}\,\E \xi^3 + 
Ct^8\, \E \xi^4,
$$
where $C$ is bounded by an absolute constant.
Summing the two expansions, we arrive at
$$
\E\,(U + V) \, = \, \frac{t^6}{6}\,\E \xi^3 + Ct^8\, \E \xi^4 + 
\ep_n(t)
$$
and therefore
$$
\int_0^1 \frac{\Delta_n(t)}{t^2}\,dt \, \leq \, \E \xi^3 + 
c\, \E \xi^4 + O(n^{-2}).
$$
Here $\E \xi^4 \geq (\E \xi^2)^2  = n^{-2}$, so the term 
$O(n^{-2})$ may be absorbed by the 4-th moment of $\xi$.
Since $\E \xi^3 \geq 0$, the bound (9.4) may be simplified to
$$
c\,\E_\theta\,\rho^2(F_\theta,F) \, \leq \,
\E \xi^3 + \E \xi^4 + 
\log n \, \int_0^{4\sqrt{\log n}} \frac{\Delta_n(t)}{t}\,dt + 
\frac{(\log n)^2}{n^2}\,m_4^4,
$$
that is,
\be
c\,\E_\theta\,\rho^2(F_\theta,F) \, \leq \,
\log n \, \int_0^{4\sqrt{\log n}} \frac{\Delta_n(t)}{t}\,dt + 
\E \xi^3 + (\log n)^2\,\E \xi^4.
\en

Turning to the remaining integral (which is most important), let us 
express it in terms of the functions $g_n(t) = J_n(t\sqrt{2n})$ and
$$
\psi(\alpha) = \int_0^T \frac{g_n(\alpha t) - g_n(t)}{t}\,dt, 
\qquad 0 \leq \alpha \leq \sqrt{2}, \ T > 1,
$$
which will be needed with $T = 4\sqrt{\log n}$ and 
$\alpha = \sqrt{1 - \xi}$. Namely, we have
\be
\int_0^T \frac{\Delta_n(t)}{t}\,dt \, = \,
\E\,\psi\big(\sqrt{1 - \xi}\big) + 
\int_0^T \frac{J_n(t\sqrt{2n}) - (J_n(t\sqrt{n}))^2}{t}\,dt.
\en

To proceed, we need to develop a Taylor expansion for 
$\xi \rightarrow \psi\big(\sqrt{1 - \xi}\big)$ around zero in powers 
of $\xi$. Recall that $g_n(t)$ represents the characteristic function
of the random variable $\sqrt{2n}\,\theta_1$ on the probability space 
$(\S^{n-1},\s_{n-1})$. This already ensures that $|g_n(t)| \leq 1$ and 
$$
|g_n'(t)| \, \leq \, \sqrt{2n}\,\E\,|\theta_1| \, \leq \, 
\sqrt{2n}\,(\E\,\theta_1^2)^{1/2} \, = \, \sqrt{2}
$$
for all $t \in \R$. Hence 
$$
|g_n(\alpha t) - g_n(t)| \leq \sqrt{2}\,|\alpha - 1|\,|t| \leq 2\,|t|,
$$ 
so that
\begin{eqnarray}
|\psi(\alpha)| 
 & \leq &
\int_0^1 \frac{|g_n(\alpha t) - g_n(t)|}{t}\,dt +
\int_1^T \frac{|g_n(\alpha t) - g_n(t)|}{t}\,dt \nonumber \\
 & \leq &
2 + 2\log T \, < \, 4\log T 
\end{eqnarray}
(since $T>e$). In addition, $\psi(1) = 0$ and
$$
\psi'(\alpha) = \int_0^T g_n'(\alpha t)\,dt = 
\frac{1}{\alpha}\,(g_n(\alpha T) - 1).
$$
Therefore, we arrive at another expression
$$
\psi(\alpha) = \int_1^\alpha \frac{g_n(Tx) - 1}{x}\,dx = 
\int_1^\alpha \frac{g_n(Tx)}{x}\,dx - \log \alpha.
$$

For $|\ep| \leq 1$, let
\bee
v(\ep) 
 & = &
\int_1^{(1 - \ep)^{1/2}} \frac{g_n(Tx)}{x}\,dx, \\
u(\ep) 
 & = & 
\psi\big((1 - \ep)^{1/2}\big) \, = \,
v(\ep) - \frac{1}{2}\,\log(1 - \ep),
\ene
so that $\E\,\psi\big(\sqrt{1 - \xi}\big) = \E\,u(\xi)$.
Applying the non-uniform bound $|g_n(t)| \leq 5\,(e^{-t^2} + e^{-n/12})$,
cf. (2.8), we have that, for $-1 \leq \ep \leq \frac{1}{2}$,
\bee
|v(\ep)| 
 & \leq &
\sup_{\frac{1}{\sqrt{2}} \leq x \leq \sqrt{2}}\,|g_n(Tx)|
\int_{\frac{1}{\sqrt{2}}}^{\sqrt{2}} \frac{1}{x}\,dx \\
 & \leq & 
\sup_{z \geq T/\sqrt{2}}\,|g_n(z)|\, \log 2 \, \leq \,
5 \log 2\,(e^{-T^2/2} + e^{-n/12}) \, \leq \, \frac{c}{n^8},
\ene
where the last inequality is specialized to the choice 
$T = 4\sqrt{\log n}$. Using the Taylor 
expansion on the same interval for the log-function, we also have
$-\log(1-\ep) \leq \ep + \frac{1}{2}\,\ep^2 + 
\frac{1}{3}\,\ep^3 + \frac{2}{3}\,\ep^4$.
Combining the two inequalities, we get
\be
u(\ep) \leq 
\frac{1}{2}\,\ep + \frac{1}{4}\,\ep^2 + \frac{1}{6}\,\ep^3 + 
\frac{1}{3}\,\ep^4 + \frac{c}{n^8}, \quad 
-1 \leq \ep \leq \frac{1}{2}.
\en

In order to involve the remaining interval 
$\frac{1}{2} \leq \ep \leq 1$ in the inequality of a similar type, recall 
that, by (10.5), $|u(\ep)| \leq 4 \log T$ for all $|\ep| \leq 1$. Hence, 
the inequality (10.6) will hold automatically for this interval, 
if we increase the coefficient 
in front of $\ep^4$ to a suitable multiple of $\log T$. As a result, 
we obtain the desired inequality on the whole segment, that is,
$$
u(\ep) \leq 
\frac{1}{2}\,\ep + \frac{1}{4}\,\ep^2 + 
\frac{1}{6}\,\ep^3 + (c\log T)\,\ep^4 + \frac{c}{n^8}.
\quad -1 \leq \ep \leq 1.
$$
In particular,
$$
\psi\big(\sqrt{1 - \xi}\big) \, \leq \,
\frac{1}{2}\,\xi + \frac{1}{4}\,\xi^2 + \frac{1}{6}\,\xi^3 + 
(c\log T)\,\xi^4 + \frac{c}{n^8},
$$
and taking the expectation, we get
\be
\E\,\psi\big(\sqrt{1 - \xi}\big) \leq 
\frac{1}{4n} + \frac{1}{6}\,\E \xi^3 + (c\log T)\, \E \xi^4,
\en
where the term $cn^{-8}$ was absorbed by the 4-th moment of $\xi$.

Now, let us turn to the integral
$$
I_n = \int_0^T \frac{J_n(t\sqrt{2n}) - (J_n(t\sqrt{n}))^2}{t}\,dt,
$$
appearing in (10.4), and recall the asymptotic formulas (10.1)-(10.2). After 
integration, the remainder term $\ep_n(t) = O\big(n^{-2} \min(1,t^4)\big)$ 
will create an error of order at most $n^{-2} \log T$, up to which $I_n$ 
is equal to
$$
-\int_0^T \frac{t^4}{2n}\,e^{-t^2}\,\frac{dt}{t} 
 \, = \,
-\frac{1}{4n}\,\Big(1 - (T^2 + 1)\,e^{-T^2}\Big) \, = \,
-\frac{1}{4n} + o(n^{-15}).
$$
Thus,
$$
I_n = -\frac{1}{4n} + O(n^{-2} \log T).
$$
Applying this expansion together with (10.7) in (10.4),
we therefore obtain that
$$
\int_0^T \frac{\Delta_n(t)}{t}\,dt \, \leq \, 
\frac{1}{6}\,\E \xi^3 + c\log T\,\E \xi^4.
$$

One can now apply this estimate in (10.3), and 
then we eventually arrive at
$$
\E_\theta\,\rho^2(F_\theta,F) \, \leq \,
c_1\, (\log n)\, \E \xi^3 + c_2\,(\log n)^2\,\E \xi^4.
$$
By (2.3) with $p = \infty$, a similar inequality remains to hold 
for the standard normal distribution function $\Phi$ in place of $F$.
This proves the inequality (1.4).
\qed


\vskip10mm
\section{{\bf Relations between $L^1$, $L^2$ and Kolmogorov 
Distances}}
\setcounter{equation}{0}

\vskip2mm
\noindent
Given a random vector $X$ in $\R^n$, let us now compare the $L^2$ 
and $L^\infty$ distances on average, between the distributions 
$F_\theta$ of the weighted sums $\left<X,\theta\right>$ and 
the typical distribution $F = \E_\theta F_\theta$. Such information 
will be needed to derive appropriate lower bounds on 
$\E_\theta\,\rho(F_\theta,F)$.

\vskip5mm
{\bf Proposition 11.1.} {\sl If $|X| \leq b\sqrt{n}$ a.s., then, for any 
$\alpha \in [1,2]$,
\be
b^{-\alpha/2}\,\E_\theta\,\omega^\alpha(F_\theta,F) \, \leq \, 
14\,(\log n)^{\alpha/4}\ \E_\theta\, \rho^\alpha(F_\theta,F) + 
\frac{8}{n^4}.
\en
}

\vskip2mm
As will be clear from the proof, at the expense of a larger coefficient 
in front of $\log n$, the last term $n^{-4}$ can be replaced by
$n^{-\beta}$ for any prescribed value of $\beta$.

A relation similar to (11.1) is also true for the Kantorovich or 
$L^1$-distance 
$$
W(F_\theta,F) = \int_{-\infty}^\infty |F_\theta(x) - F(x)|\,dx
$$
in place of $L^2$. We state it for the case $\alpha = 1$.

\vskip5mm
{\bf Proposition 11.2.} {\sl If $|X| \leq b\sqrt{n}$ a.s., then
\be
\E_\theta\,W(F_\theta,F) \, \leq \, 
14\,b\sqrt{\log n}\ \E_\theta\, \rho(F_\theta,F) + \frac{8b}{n^4}.
\en
}

\vskip2mm
{\bf Proof.} Put $R_\theta(x) = F_\theta(-x) + (1-F_\theta(x))$ 
for $x > 0$ and define similarly $R$ on the basis of $F$. Using 
\bee
(F_\theta(-x) - F(-x))^2 
 & \leq &
F_\theta(-x)^2 + F(-x)^2,  \\
(F_\theta(x) - F(x))^2 
 & \leq &
(1-F_\theta(x))^2 + (1-F(x))^2,
\ene
we have
$$
(F_\theta(-x) - F(-x))^2 + (F_\theta(x) - F(x))^2 \leq 
R_\theta(x)^2 + R(x)^2.
$$
Hence, given $T > 0$ (to be specified later on), we have
\bee
\omega^2(F_\theta,F)
 & = &
\int_{-T}^T (F_\theta(x) - F(x))^2\,dx + 
\int_{|x| \geq T} (F_\theta(x) - F(x))^2\,dx \\ 
 & \leq &
2T \rho^2(F_\theta,F) + \int_T^\infty R_\theta(x)^2\,dx  + 
\int_T^\infty R(x)^2\,dx.
\ene
It follows that, for any $\alpha \in [1,2]$,
$$
\omega^\alpha(F_\theta,F) \leq (2T)^{\frac{\alpha}{2}}\, 
\rho^\alpha(F_\theta,F) + 
\Big(\int_T^\infty R_\theta(x)^2\,dx\Big)^{\frac{\alpha}{2}}
 + \Big(\int_T^\infty R(x)^2\,dx\Big)^{\frac{\alpha}{2}}
$$
and therefore, by Jensen's inequality,
\bee
\E_\theta\,\omega^\alpha(F_\theta,F) 
 & \leq &
(2T)^{\frac{\alpha}{2}}\ \E_\theta\, \rho^\alpha(F_\theta,F) \\
 & & + \ 
\Big(\int_T^\infty \E_\theta\, 
R_\theta(x)^2\,dx\Big)^{\frac{\alpha}{2}} + 
\Big(\int_T^\infty R(x)^2\,dx\Big)^{\frac{\alpha}{2}}.
\ene

Next, by Markov's inequality, for any $x > 0$ and $p \geq 1$,
$$
R_\theta(x)^2 \leq 
\Big(\frac{\E\,|\left<X,\theta\right>|^p}{x^p}\Big)^2 \leq 
\frac{\E\,|\left<X,\theta\right>|^{2p}}{x^{2p}}
$$
and
$$
\E_\theta R_\theta(x)^2 \leq 
\Big(\frac{\E\,|\left<X,\theta\right>|^p}{x^p}\Big)^2 \leq 
\frac{\E_\theta\, \E\,|\left<X,\theta\right>|^{2p}}{x^{2p}}.
$$
Since $R = \E_\theta R_\theta$, a similar inequality holds 
true for $R$ as well (by Cauchy's inequality). Hence
$$
\E_\theta\,\omega^\alpha(F_\theta,F) \, \leq \,
(2T)^{\frac{\alpha}{2}}\ \E_\theta\, \rho^\alpha(F_\theta,F) + 
2\,\Big(\E_\theta\,\E\,|\left<X,\theta\right>|^{2p} 
\int_T^\infty
\frac{1}{x^{2p}}\,dx\Big)^{\frac{\alpha}{2}}.
$$
When $\theta = (\theta_1,\dots,\theta_n)$ is treated 
as a random vector with distribution $\s_{n-1}$, which is 
independent of $X$, the inner product $\left<X,\theta\right>$ 
has the same distribution as the random variable 
$|X|\,\theta_1$. Therefore, recalling Lemma 2.5 and using
the assumption $|X| \leq b\sqrt{n}$ a.e., we have
$$
\E_\theta\,\E\,|\left<X,\theta\right>|^{2p} \, = \, 
\E\,|X|^{2p} \ \E_\theta\, |\theta_1|^{2p} \, \leq \,
2\,(2b^2 p)^p,
$$
so that
$$
2\,\Big(\E_\theta\, \int_T^\infty
\frac{\E\,
|\left<X,\theta\right>|^{2p}}{x^{2p}}\,dx\Big)^{\frac{\alpha}{2}}
 \, \leq \, 
 \frac{2^{\frac{\alpha}{2} + 1}}{(2p-1)^{\frac{\alpha}{2}}}\,
\frac{(2b^2 p)^{\frac{\alpha p}{2}}}{T^{\frac{\alpha (2p-1)}{2}}}.
$$

Thus,
$$
\E_\theta\,\omega^\alpha(F_\theta,F) \, \leq \,
(2T)^{\frac{\alpha}{2}}\ \E_\theta\, \rho^\alpha(F_\theta,F) + 
\frac{2^{\frac{\alpha}{2} + 1}}{(2p-1)^{\frac{\alpha}{2}}}\ 
T^{\frac{\alpha}{2}}\, 
\Big(\frac{2b^2 p}{T^2}\Big)^{\frac{\alpha p}{2}}.
$$
Let us choose $T = 2b\sqrt{p}$ in which case the above inequality 
becomes
$$
\E_\theta\,\omega^\alpha(F_\theta,F) \, \leq \,
(4b\sqrt{p})^{\frac{\alpha}{2}}\ \E_\theta\, \rho^\alpha(F_\theta,F) + 
\frac{2^{\alpha + 1}}{(2p-1)^{\frac{\alpha}{2}}}\ 
(b\sqrt{p})^{\frac{\alpha}{2}}\ 2^{-\frac{\alpha p}{2}}.
$$
To simplify, one can use $\sqrt{p} \leq 2p-1$ 
for $p \geq 1$ together with $2^{\alpha + 1} \leq 8$ 
and $2^{-\frac{\alpha p}{2}} \leq 2^{-\frac{p}{2}}$
(since $1 \leq \alpha \leq 2$), which leads to
$$
\E_\theta\,\omega^\alpha(F_\theta,F) \, \leq \,
(4b\sqrt{p})^{\frac{\alpha}{2}}\ \E_\theta\, \rho^\alpha(F_\theta,F) + 
8\, b^{\frac{\alpha}{2}}\ 2^{-p/2}.
$$
Finally, choosing $p = p_n = (8\,\log n)/\log 2$, we arrive at (11.1).

Now, turning to (11.2), we use the  same functions $R_\theta$ and 
$R$ as before and write
\bee
W(F_\theta,F)
 & = &
\int_{-T}^T |F_\theta(x) - F(x)|\,dx + 
\int_{|x| \geq T} |F_\theta(x) - F(x)|\,dx \\ 
 & \leq &
2T \rho(F_\theta,F) + \int_T^\infty R_\theta(x)\,dx  + 
\int_T^\infty R(x)\,dx,
\ene
which gives
$$
\E_\theta\,W(F_\theta,F) \, \leq \, 
2T\, \E_\theta\, \rho(F_\theta,F) + 2 \int_T^\infty R(x)\,dx.
$$
By Markov's inequality, for any $x > 0$ and $p > 1$,
$$
R_\theta(x) \leq \frac{\E\,|\left<X,\theta\right>|^p}{x^p}, 
\quad R(x) = \E_\theta R_\theta(x) \leq 
\frac{\E_\theta\, \E\,|\left<X,\theta\right>|^p}{x^p}.
$$
Hence
$$
\E_\theta\,W(F_\theta,F) \, \leq \, 
2T\, \E_\theta\, \rho(F_\theta,F) + 
2\,\E_\theta\, \E\,|\left<X,\theta\right>|^p \int_T^\infty
\frac{1}{x^p}\,dx.
$$
Here, one may use once more the bound (2.10), which yields
$$
\E_\theta\,\E\,|\left<X,\theta\right>|^p \, = \, 
\E\,|X|^p \ \E_\theta\, |\theta_1|^p \, \leq \,
2\,\big(b^2 p\big)^{p/2}
$$
and
$$
\E_\theta\,W(F_\theta,F) \, \leq \, 
2T\, \E_\theta\, \rho(F_\theta,F) + \frac{4}{p-1}\,
\frac{(b^2 p)^{p/2}}{T^{p-1}}.
$$
Let us take $T = 2b\sqrt{p}$ in which case the above inequality 
becomes
$$
\E_\theta\,W(F_\theta,F) \, \leq \, 
4b\sqrt{p}\ \E_\theta\, \rho(F_\theta,F) + 
8b\,\frac{\sqrt{p}}{p-1}\ 2^{-p}.
$$
Here we arrive at (11.2), by choosing again $p = p_n$ and 
using $\sqrt{p_n} < p_n - 1$.
\qed


\vskip10mm
\section{{\bf Lower Bounds. Proof of Theorem 1.3}}
\setcounter{equation}{0}

\vskip2mm
\noindent
A lower bound on $\E_\theta\,\rho^2(F_\theta,\Phi)$ which would be close 
to the upper bound (1.4) may be given with the help of the lower bound 
on $\E_\theta\, \omega^2(F_\theta,\Phi)$. More precisely, this 
can be done in the case where the quantity 
$\frac{1}{n^{3/2}}\, m_3^3 + \frac{1}{n^2}\, m_4^4$
asymptotically dominates $n^{-2}$ (in particular, when $m_4$
is essentially larger than 1). 
Combining the asymptotic expansion (1.3) of Theorem 1.1 with 
the bound (11.1) of Proposition 11.1 for $\alpha=2$ and $b=1$, and
recalling the second relation in (2.3) on the normal approximation 
for the typical distribution $F$, we therefore obtain:

\vskip5mm
{\bf Proposition 12.1.} {\sl If $X$ is an isotropic random vector 
in $\R^n$ with mean zero and such that $|X| = \sqrt{n}$ a.s., then
\be
\sqrt{\log n} \ \E_\theta\,\rho^2(F_\theta,\Phi) \, \geq \,
\frac{c_1}{n^{3/2}}\, m_3^3 + 
\frac{c_2}{n^2}\, m_4^4 - \frac{c_3}{n^2}.
\en
}

\vskip2mm
The relation (11.2) for the Kantorovich distance $W$ may 
be used to answer the following question: Is it possible 
to sharpen the lower bound (12.1) by replacing 
$\E_\theta\,\rho^2(F_\theta,\Phi)$ with $\E_\theta\,\rho(F_\theta,\Phi)$?
To this aim, we will need an additional information about 
moments of $\omega(F_\theta,F)$ of order higher than 2.

\vskip5mm
{\bf Lemma 12.2.} {\sl If $X$ is isotropic and satisfies 
$|X| \leq b\sqrt{n}$, then 
\be
c\,\left(\E_\theta\,\omega^3(F_\theta,F)\right)^{1/3} \, \leq \, 
(1 + \sigma_4)\sqrt{b}\ \frac{(\log n)^{5/4}}{\sqrt{n}}.
\en
}

\vskip2mm
{\bf Proof.} For any distribution function $G$ with finite first 
absolute moment, the function on the unit sphere $\S^{n-1}$ of the
form $g(\theta) = W(F_\theta,G)$ has a Lipschitz semi-norm
$\|g\|_{\rm Lip} \leq 1$. Therefore, it admits a subgaussian large 
deviation bound
\be
\s_{n-1}\big\{W(F_\theta,G) \geq m + r\big\} \, \leq \, e^{-(n-1)r^2/2},
\qquad r \geq 0,
\en
where $m = \E_\theta\, W(F_\theta,G)$. Indeed, consider 
the elementary representation
\bee
W(F_\theta,G) 
 & \equiv &
\int_{-\infty}^\infty |F_\theta(x) - G(x)|\,dx \\ 
 & = &
\sup_u \bigg[\int_{-\infty}^\infty u\,dF_\theta -
\int_{-\infty}^\infty u\,dG\bigg],
\ene
where the supremum is running over all functions $u$ on $\R$
with $\|u\|_{\rm Lip} \leq 1$. For any such $u$,
$$
H_u(\theta) = \int_{-\infty}^\infty u\,dF_\theta \, = \, 
\E\,u(\left<X,\theta\right>)
$$
is Lipschitz on $\R^n$ and therefore on $\S^{n-1}$. Moreover,
$\|g\|_{\rm Lip} \leq \sup_u \|H_u\|_{\rm Lip} \leq 1$.

Hence, (12.3) is fulfilled as a consequence of fact
that the logarithmic Sobolev constant for the uniform distribution on 
the unit sphere is equal to $n-1$ (cf. \cite{L2}). In particular,
for any $r \geq 0$,
$$
\s_{n-1}\big\{W(F_\theta,F) \geq m + r\big\} \, \leq \, e^{-(n-1)r^2/2}
$$
with $m = \E_\theta\, W(F_\theta,F)$. In turn, the latter ensures 
that, for any $p \geq 2$,
\be
\Big(\E_\theta\, W(F_\theta,F)^p\Big)^{1/p} \leq 
m + \frac{\sqrt{p}}{\sqrt{n-1}}.
\en
For the proof, put $\xi = (W(F_\theta,F) - m)^+$. Using 
$\Gamma(x+1) \leq x^x$ with $x = p/2 \geq 1$, we have
\bee
\E_\theta\, \xi^p 
 & = &
\int_0^\infty \s_{n-1}\{\xi \geq r\}\,dr^p
 \, \leq \,
\int_0^\infty e^{-(n-1)\,r^2/2}\,dr^p \\
 & = & 
\Big(\frac{\sqrt{2}}{\sqrt{n-1}}\Big)^p\ 
\Gamma\Big(\frac{p}{2} + 1\Big)
 \, \leq \,
\bigg(\frac{\sqrt{p}}{\sqrt{n-1}}\,\bigg)^p \, \equiv \, A^p \quad
(A \geq 0).
\ene
Thus, $\|\xi\|_p = (\E_\theta\, \xi^p)^{1/p} \leq A$.
Since $W(F_\theta,F) \leq \xi + m$, we conclude, by the
triangle inequality, that
$$
\|W(F_\theta,F)\|_p \leq \|\xi\|_p + m \leq A + m,
$$ 
that is, (12.4) holds.

Let us proceed with one elementary general inequality, 
connecting the three distances,
\bee
\omega^2(F_\theta,F)
 & = &
\int_{-\infty}^\infty (F_\theta(x) - F(x))^2\,dx \\
 & \leq &
\int_{-\infty}^\infty \sup_x |F_\theta(x) - F(x)| \ 
|F_\theta(x) - F(x)|\,dx  \, = \, \rho(F_\theta,F)\, W(F_\theta,F).
\ene
Putting $\omega = \omega(F_\theta,F)$, $W = W(F_\theta,F)$, 
$\rho = \rho(F_\theta,F)$, we thus have 
$
\omega^3 \leq W^{3/2} \rho^{3/2}
$
and, by H\"older's inequality with exponents $p = 4$ and $q = 4/3$,
$$
\|\omega\|_3 = \big(\E_\theta\,\omega^3\big)^{1/3} \leq
\big(\E_\theta\,W^6\big)^{1/12} \, 
\big(\E_\theta\,\rho^2\big)^{1/4}.
$$
By (12.4) with $p=6$, we have 
$$
\big(\E_\theta\,W^6\big)^{1/6} \leq 
\E_\theta\,W + \frac{4}{\sqrt{n}},
$$ 
so that
$$
\|\omega\|_3 \leq \Big(\E_\theta\,W + 
\frac{4}{\sqrt{n}}\,\Big)^{1/2} \, 
\big(\E_\theta\,\rho^2\big)^{1/4}.
$$
Applying Proposition 11.2 and noting  that necessarily $b \geq 1$
in the isotrpic case, we get
$$
\|\omega\|_3 \leq
4\sqrt{b}\,\Big(\sqrt{\log n}\ \E_\theta\, \rho +
\frac{1}{\sqrt{n}}\,\Big)^{1/2} \, \big(\E_\theta\,\rho^2\big)^{1/4}.
$$
Here we employ the inequality (9.3) with $F$ in place of 
$\Phi$, i.e.
$$
\E_\theta\,\rho(F_\theta,F) \leq
\big(\E_\theta\,\rho^2(F_\theta,F)\big)^{1/2} \leq
c\,(1 + \sigma_4)\, \frac{\log n}{\sqrt{n}}.
$$
Since the last expression dominates the term $\frac{1}{\sqrt{n}}$,
it follows that
$$
\|\omega\|_3 \leq c\sqrt{b}\,
\Big(\sqrt{\log n}\ (1 + \sigma_4)\, \frac{\log n}{\sqrt{n}}\,\Big)^{1/2} \, 
\Big((1 + \sigma_4)\, \frac{\log n}{\sqrt{n}}\Big)^{1/2},
$$
and we arrive at the upper bound (12.2).
\qed

\vskip5mm
Let us now explain how this bound can be used to refine the lower 
bound (12.1). The argument is based on the following general 
elementary observation. Given a random variable $\xi$, introduce 
the $L^p$-norms $\|\xi\|_p = (\E\,|\xi|^p)^{1/p}$.

\vskip5mm
{\bf Lemma 12.3.} {\sl If $\xi \geq 0$ 
with $0 < \|\xi\|_3 < \infty$, then
\be
\E\,\xi\geq \frac{1}{\E\,\xi^3}\,(\E\,\xi^2)^2.
\en
Moreover,
\be
\P\Big\{\xi \geq \frac{1}{\sqrt{2}}\,\|\xi\|_2\Big\} \, \geq \, 
\frac{1}{8}\,\Big(\frac{\|\xi\|_2}{\|\xi\|_3}\Big)^6.
\en
}

\vskip2mm
Thus, in the case where $\|\xi\|_2$ and $\|\xi\|_3$ are 
equivalent within not too large factors, $\|\xi\|_1$ will be of a similar 
order. Moreover, $\xi$ cannot be much smaller than its mean 
$\E \xi$ on a large part of the probability space (where it was defined).

\vskip5mm
{\bf Proof.} Let $\xi$ be defined on the probability space
$(\Omega,{\mathfrak F},\P)$. By homogeneity with respect to $\xi$,
we may assume that
$\E \xi = 1$, so that $dQ = \xi d\P$ is a probability measure. Then, 
(12.5) follows from the Cauchy inequality $(\E_Q \xi)^2 \leq \E_Q \xi^2$
on the space $(\Omega,{\mathfrak F},Q)$.

To prove (12.6), given $r > 0$, let $p = \P\{\xi \geq r\}$. 
By H\"older's inequality with exponents $3/2$ and $3$,
$$
\E\,\xi^2\,1_{\{\xi \geq r\}} \leq 
\big(\E\,\xi^3\big)^{2/3}\, p^{1/3}.
$$
Hence, choosing $r = \frac{1}{\sqrt{2}}\,\|\xi\|_2$, we get
\bee
\E\,\xi^2 
 & = &
\E\,\xi^2\,1_{\{\xi \geq r\}} + \E\,\xi^2\,1_{\{\xi < r\}} \\
 & \leq &
\big(\E\,\xi^3\big)^{2/3}\, p^{1/3} + r^2 \, = \, 
\big(\E\,\xi^3\big)^{2/3}\, p^{1/3} + \frac{1}{2}\,\E\,\xi^2.
\ene
Hence $p^{1/3} \geq \frac{1}{2\,(\E\,\xi^3)^{2/3}}\,\E\,\xi^2$
which is the desired bound (12.6).
\qed

\vskip5mm
We now combine Lemma 12.2 with Lemma 12.3 which is applied 
on the unit sphere to $\xi(\theta) = \omega(F_\theta,F)$ viewed 
as a random variable on the probability space $(\S^{n-1},\s_{n-1})$. 
Recall that $b \geq 1$ in the isotropic case.

\vskip5mm
{\bf Proposition 12.4.} {\sl Let $X$ be an isotropic random vector 
in $\R^n$ such that $|X| \leq b\sqrt{n}$ a.s. Assume that
$$
\E_\theta\,\omega^2(F_\theta,F) \, \geq \, 
\frac{D}{n}
$$
with some $D>0$. Then 
\be
\E_\theta\,\omega(F_\theta,F)\,\geq\, 
\frac{c}{(1 + \sigma_4)^3\, b^{\frac{3}{2}}}\,
\frac{D^2}{(\log n)^{\frac{15}{4}} \sqrt{n}}.
\en
Moreover, 
$$
\s_{n-1}\Big\{\omega(F_\theta,F) \geq \frac{1}{\sqrt{2n}}\sqrt{D}\Big\}
 \, \geq \, \frac{c}{(1 + \sigma_4)^6\, b^3}\,
\frac{D^3}{(\log n)^{\frac{15}{2}}}.
$$
}

\vskip5mm
{\bf Proof of Theorem 1.3.}
The lower bound (12.7) implies a similar assertion about 
the Kolmogorov distance. Indeed, by Proposition 11.1 with 
$\alpha = 1$, we have
$$
\frac{1}{\sqrt{b}}\, \E_\theta\,\omega(F_\theta,F) 
 \, \leq \, 14\,(\log n)^{1/4}\ 
 \E_\theta\, \rho(F_\theta,F) + \frac{8}{n^4}.
$$
Using $\frac{8}{n^4} < \frac{1}{n^3} \cdot 14\,(\log n)^{1/4}$, 
we therefore obtain that
\bee
\E_\theta\,\rho(F_\theta,F)
 & \geq &
\frac{1}{14 \sqrt{b}\ (\log n)^{1/4}}\ \E_\theta\,\omega(F_\theta,F) -
\frac{1}{n^3} \\
 & \geq &
\frac{c}{(1 + \sigma_4)^3\,b^2}\,
\frac{D^2}{(\log n)^4\, \sqrt{n}} - \frac{1}{n^3}.
\ene
To replace $F$ with $\Phi$, it remains to recall the bound 
$\rho(F,\Phi) \leq \frac{c}{n}\,(1 + \sigma_4^2)$, cf. (2.3).
\qed

\vskip5mm
In the isotropic case with $|X|^2 = n$ a.s., the above lower bound is further 
simplified to
$$
\E_\theta\,\rho(F_\theta,F)\,\geq\, 
\frac{cD^2}{(\log n)^4\,\sqrt{n}} - \frac{1}{n^3}.
$$
On the other hand, let us
note that the rates for the normal approximation of $F_\theta$ 
that are better than $1/n$ (on average) cannot be obtained 
under the support assumption as above. That is, if $|X| = \sqrt{n}$ a.s., 
then
$$
\E_\theta\,\rho(F_\theta,\Phi)\,\geq\, \frac{c}{n}.
$$
Indeed, using the convexity of the distance function 
$G \rightarrow \rho(G,\Phi)$ and applying Jensen's inequality, 
we have that $\E_\theta\,\rho(F_\theta,\Phi) \geq \rho(F,\Phi)$.
It remains to appeal to Proposition 2.6.


\vskip10mm
\section{{\bf Functional Examples}}
\setcounter{equation}{0}

\vskip2mm
\noindent
{\bf 13.1.} For the trigonometric system as in item (i) of the Introduction
(with $n$ even), the linear forms
$$
\left<X,\theta\right> =  \sqrt{2}\, \sum_{k=1}^{\frac{n}{2}} 
\big(\theta_{2k-1} \cos(kt) + \theta_{2k} \sin(kt)\big),
\quad \theta = (\theta_1,\dots,\theta_n) \in \S^{n-1},
$$
represent trigonometric polynomials of degree at most $\frac{n}{2}$.
The normalization $\sqrt{2}$ is chosen in order to meet the
requirement that the random vector $X$ is isotropic with respect
to the normalized Lebesgue measure $\P$ on $\Omega = (-\pi,\pi)$.
Moreover, in this case $|X| = \sqrt{n}$, so that $\sigma_4 = 0$. 
Hence, by Theorem 1.1, we have the upper bounds (1.6). 
On the other hand, since for all $k \leq \frac{n}{2}$
$$
|X_k(t) - X_k(s)| \leq k\sqrt{2}\ |t - s| \leq 
\frac{n}{\sqrt{2}}\,|t - s|, \qquad t,s \in \Omega,
$$
the Lipschitz condition (7.1) is fulfilled with $L(t) =  \frac{t}{\sqrt{2}}$.
Hence, Proposition 7.1 is applicable and yields the lower bound 
$$
\E_\theta\, \omega^2(F_\theta,\Phi) \,\geq\, 
\frac{c_1}{n} - \frac{c_2}{n^2} \,\geq\, \frac{c_3}{n},
$$
where in the last inequality we assume that $n \geq n_0$ for some 
universal integer $n_0$. This restriction may be dropped, since the
distances $\omega^2(F_\theta,\Phi)$ are bounded away from zero
for $n < n_0$ uniformly over all $\theta \in \S^{n-1}$, just due 
to the property that the distributions $F_\theta$ are supported on 
the bounded interval $[-\sqrt{n_0},\sqrt{n_0}]$. Note that 
the above lower estimate (may also be 
obtained by applying Theorem 1.1. Thus, for all $n \geq 2$,
\be
\frac{c_0}{n}\,\leq\, \E_\theta\, \omega^2(F_\theta,\Phi) 
 \,\leq\, \frac{c_1}{n}.
\en

Applying Proposition 12.4, we obtain similar
bounds for the $L^1$-norm (modulo logarithmic factors). 
Namely, it gives
\be
\frac{c_0}{(\log n)^\frac{15}{4}\sqrt{n}} \leq 
\E_\theta\, \omega(F_\theta,\Phi) \leq \frac{c_1}{\sqrt{n}}.
\en
We also get an analogous pointwise lower bound on the 
``essential" part of the unit sphere. 

A similar statement is also true for the Kolmogorov distance.
Here, the upper bound is provided in Proposition 9.1, while
the lower bound is obtained when combining
Theorem 1.3 with the left inequality in (13.1). That is,
\be
\frac{c_0}{(\log n)^4 \sqrt{n}} \leq 
\E_\theta\, \rho(F_\theta,\Phi)  \leq 
\big(\E_\theta\, \rho^2(F_\theta,\Phi)\big)^{1/2} \leq \frac{c_1\log n}{\sqrt{n}}.
\en

{\bf 13.2.} 
Analogous results remain true for the cosine trigonometric system
$X = (X_1,\dots,X_n)$ as in item (ii). Due to the normalization $\sqrt{2}$, 
the distribution of $X$ is isotropic in $\R^n$. The property 
$|X| = \sqrt{n}$ is not true anymore; however, there is a pointwise bound 
$|X| \leq \sqrt{2n}$.  In addition, the variance functional $\sigma_4^2$
does not depend on $n$. Indeed, write
$$
X_k^2 = 2\cos^2(kt) = 1 + \cos(2kt) = 
1 + \frac{e^{2ikt} + e^{-2ikt}}{2},
$$
so that
$$
2\,(|X|^2 - n) \ = \sum_{0 < |k| \leq n} e^{2ikt}, \qquad
4\,(|X|^2 - n)^2 \ = \sum_{0 < |k|, |l| \leq n} e^{2i(k+l)t}.
$$
It follows that
$$
4\,\Var(|X|^2) \ = \sum_{0 < |k|, |l| \leq n} \E\,e^{2i(k+l)t}
 \ = \sum_{0< |k| \leq n,\, l = -k} 1 \, = \, 2n.
$$
Hence
$$
\sigma_4^2 = \frac{1}{n}\,\Var(|X|^2) = \frac{1}{2}.
$$

As before, the Lipschitz condition is fulfilled with the function 
$L(t) =  t\sqrt{2}$. Therefore, with similar arguments we
obtain all the bounds (13.1)-(13.3).

Let us also note that the sums $\sum_{k=1}^n \cos(kt)$ 
remain bounded for growing $n$ (for any fixed $0 < t < \pi$). 
Hence the normalized sums
$$
S_n = \frac{1}{\sqrt{n}} \sum_{k=1}^n X_k = 
\frac{\sqrt{2}}{\sqrt{n}}\,\sum_{k=1}^n \cos(kt),
$$
which correspond to $\left<X,\theta\right>$ with equal 
coefficients, are convergent to zero pointwise on $\Omega$ 
as $n \rightarrow \infty$. In particular, they fail to satisfy 
the central limit theorem.

{\bf 13.3.} 
An example closely related to the cosine trigonometric system
is represented by the normalized Chebyshev's polynomials $X_k$
as in item (iii), which we consider for $k = 1,2,\dots,n$. These 
polynomials are orthonormal on the interval $\Omega = (-1,1)$ 
with respect to the probability measure 
$$
\frac{d\P(t)}{dt} = \frac{1}{\pi \sqrt{1 - t^2}}, \quad -1<t<1,
$$ 
cf. e.g. \cite{K-Sa}. Similarly to 13.2, for 
the random vector $X = (X_1,\dots,X_n)$ we find that
$$
4\,(|X|^2 - n)^2 \ = \sum_{0 < |k|, |l| \leq n} 
\exp\{2i(k+l)\, \arccos t\}.
$$
It follows that
$$
4\,\Var(|X|^2) \ = 
\sum_{0 < |k|, |l| \leq n} \E\,\exp\{2i(k+l)\,\arccos t\}
 \ = \sum_{0< |k| \leq n} 1 \, = \, 2n,
$$
so that $\sigma_4^2 = \frac{1}{n}\,\Var(|X|^2) = \frac{1}{2}$.
In addition, for all $k \leq n$,
$$
|X_k(t) - X_k(s)| \leq k\sqrt{2}\ |\arccos t - \arccos s|, 
 \qquad t,s \in \Omega,
$$
which implies that the Lipschitz condition is fulfilled with the function 
$L(t) =  \sqrt{2}\,\arccos t$. As a result, 
we obtain the bounds (13.1)-(13.3) as well.

{\bf 13.4.} Turning to item (iv), consider the functions of the form 
$$
X_k(t,s) = \Psi(kt + s),
$$ 
assuming that $\Psi$ is a 1-periodic measurable function on the real line 
such that 
$$
\int_0^1 \Psi(x)\,dx = 0  \quad {\rm and} \quad 
\int_0^1 \Psi(x)^2\,dx = 1.
$$
These conditions ensure that the random vector 
$X = (X_1,\dots,X_n)$ is isotropic in $\R^n$ with respect to 
the Lebesgue measure $\P$ on the square 
$\Omega = (0,1) \times (0,1)$, with $\E X_k = 0$. In fact, as was 
emphasized in \cite{B-G}, $\{X_k\}_{k=1}^\infty$ represents a strictly 
stationary sequence of pairwise independent random variables on $\Omega$. 
The latter implies in particular that, if $\Psi$ has finite 4-th moment 
on $(0,1)$, the variance functional
$$
\sigma_4^2 = \frac{1}{n}\,\Var(|X|^2) = 
\int_0^1 \Psi(x)^4\,dx - 1
$$
is finite and does not dependent on $n$. Hence, by Theorem 1.1, 
cf. (1.6), the upper bounds in (13.1)-(13.3) hold true
with a constant $c_1$ depending on the 4-th moment of $\Psi$
on $(0,1)$.

In addition, if the function $\Psi$ has finite Lipschitz constant
$\|\Psi\|_{\rm Lip}$, then for all $(t_1,t_2)$ and $(s_1,s_2)$ in 
$\Omega$,
$$
|X_k(t_1,t_2) - X_k(s_1,s_2)| \, \leq \, \|\Psi\|_{\rm Lip}\,
\big(k\,|t_1 - s_1| + |t_2 - s_2|\big).
$$
This means that the  Lipschitz condition (7.5) is 
fulfilled with linear functions $L_1$ and $L_2$. Hence, 
one may apply Proposition 7.5 giving the lower bound
$$
\E_\theta\, \omega^2(F_\theta,F) \,\geq\, 
\frac{c_\Psi}{n} - \frac{c\,(1 + \sigma_4^4)}{n^2}
$$
in full analogy with item (i). Hence 
$\E_\theta\, \omega^2(F_\theta,\Phi) \,\geq\, \frac{c_\Psi'}{n}$
for all $n \geq n_0$, where the positive constants $c_\Psi$, $c_\Psi'$, 
and an integer $n_0 \geq 1$ depend on the distribution 
of $\Psi$ only. Since the collection $\{F_\theta\}$ is separated
from $\Phi$ in the weak sense for $n<n_0$ (by the uniform 
boundedness of $X_k$'s), the latter bound holds true for 
all $n \geq 2$. Also, as Lipschitz functions on $(0,1)$ are bounded, 
we have $|X| \leq b\sqrt{n}$ with $b = \sup_x |f(x)|$, and one may
apply Theorem 1.3. 

\vskip2mm
Let us summarize: {\sl The upper bounds in $(13.1)-(13.3)$ hold true,
if $\Psi$ has finite 4-th moment under the uniform distribution
on $(0,1)$. The lower bounds hold under an additional
assumption that $\Psi$ has a finite Lipschitz semi-norm
(with constants depending on $\Psi$ only).
}

\vskip2mm
Choosing, for example, $\Psi(t) = \cos t$, we obtain the system 
$X_k(t,s) = \cos(kt + s)$, which is closely related to the cosine 
trigonometric system. The main difference is however the property 
that $X_k$'s are now pairwise independent. Nevertheless, 
the normalized sums $\frac{1}{\sqrt{n}} \sum_{k=1}^n \cos(kt + s)$
fail to satisfy the central limit theorem.


\vskip10mm
\section{{\bf The Walsh System; Empirical Measures}}
\setcounter{equation}{0}

\vskip2mm
\noindent
{\bf 14.1.} The Walsh system on the discrete cube $\Omega = \{-1,1\}^d$ 
with the uniform counting measure $\P$ as in item (v) in Introduction 
forms a complete orthonormal system in $L^2(\Omega,\P)$. Note that each $X_\tau$ 
with $\tau \neq \emptyset$ is a symmetric Bernoulli random variable 
taking the values $-1$ and $1$ with probability $\frac{1}{2}$. 
For simplicity, we exclude from this family the constant
$X_{\emptyset} = 1$ and consider $X = \{X_\tau\}_{\tau \neq \emptyset}$ 
as a random vector in $\R^n$ of dimension $n = 2^d - 1$. As before,
$F_\theta$ denotes the distribution function of the linear form
$$
\left<X,\theta\right> = 
\sum_{\tau \neq \emptyset} \theta_\tau X_\tau, 
\quad \theta = 
\{\theta_\tau\}_{\tau \neq \emptyset} \in \S^{n-1}.
$$

Since $|X_\tau| = 1$ and thus $|X| = \sqrt{n}$, for the 
study of the asymptotic behavior of the $L^2$-distance 
$\omega(F_\theta,\Phi)$ on average, one may apply Theorem 1.1.
Let $Y$ be an independent copy of $X$, which we realize 
on the product space $\Omega^2 = \Omega \times \Omega$ 
with product measure $\P^2 = \P \times \P$ by 
$$
X_\tau(t,s) = \prod_{k \in \tau} t_k, \ \ 
Y_\tau(t,s) = \prod_{k \in \tau} s_k
\qquad 
t = (t_1,\dots,t_d), \ s = (s_1,\dots,s_d) \in \Omega.
$$
Then the inner product
$$
\left<X,Y\right> \, = \, \sum_{\tau \neq \emptyset} 
X_\tau(t,s) Y_\tau(t,s) \, = \,
-1 + \prod_{k = 1}^d \, (1 + t_k s_k)
$$
takes only two values, namely $2^d - 1$ in the case $t=s$, 
and $-1$ if $t \neq s$. Hence
$$
\E \left<X,Y\right>^3 \, = \, 
(2^d - 1)^3\, 2^{-d} + (1 - 2^{-d}) \, = \, 
\frac{n^3}{n+1} + \Big(1 - \frac{1}{n+1}\Big) \, \sim \, n^2
$$
and
$$
\E \left<X,Y\right>^4 \, = \, (2^d - 1)^4\, 2^{-d} + (1 - 2^{-d}) 
 \, = \, \frac{n^4}{n+1} + \Big(1 - \frac{1}{n+1}\Big) \, \sim \, n^3.
$$
In other words, $m_3^3 \sim \sqrt{n}$ and $m_4^4 \sim n$ as
$n \rightarrow \infty$. As a result, we may conclude that all 
inequalities in (13.1)-(13.3) are fulfilled for this system as well.

{\bf 14.2.} 
Here is another interesting example leading to the similar rate 
of normal approximation. Let $e_1,\dots,e_n$ denote the canonical 
basis in $\R^n$. Assuming that the random vector 
$X = (X_1,\dots,X_n)$ takes only $n$ values, 
$\sqrt{n}\,e_1, \dots, \sqrt{n}\,e_n$,
each with probability $1/n$, the linear form 
$\left<X,\theta\right>$ also takes $n$ values, namely, 
$\sqrt{n}\,\theta_1, \dots, \sqrt{n}\,\theta_n$, each 
with probability $1/n$, for any 
$\theta = (\theta_1,\dots,\theta_n) \in \S^{n-1}$.
That is, as a measure, the distribution of 
$\left<X,\theta\right>$ is described as
$$
F_\theta = 
\frac{1}{n}\,\sum_{k=1}^n \delta_{\sqrt{n}\,\theta_k},
$$
which may be viewed as an empirical measure based on the observations 
$Z_k = \sqrt{n}\,\theta_k$, $k = 1,\dots,n$. Each $Z_k$ is almost 
standard normal, while jointly they are nearly independent (we have 
already considered in detail its characteristic functions $J_n(t\sqrt{n})$).

Just taking a short break, let us recall that when $Z_k$ are 
indeed standard normal and independent, it is well-known that 
the empirical measures
$G_n = \frac{1}{n}\,\sum_{k=1}^n \delta_{Z_k}$
approximate the standard normal law $\Phi$ with rate 
$1/\sqrt{n}$ with respect to the Kolmogorov distance. More precisely, 
$\E\, G_n = \Phi$ and there is a subgaussian deviation bound
(cf. \cite{Mas})
$$
\P\big\{\sqrt{n}\,\rho(G_n,\Phi) \geq r\big\} \leq 
2 e^{-2r^2}, \qquad r \geq 0.
$$
In particular, $\E\,\rho(G_n,\Phi) \leq \frac{c}{\sqrt{n}}$.
Note that the characteristic function
$g_n(t) = \frac{1}{n}\,\sum_{k=1}^n e^{itZ_k}$ 
of the measure $G_n$ has mean $g(t) = e^{-t^2/2}$ and variance 
$$
\E\,|g_n(t) - g(t)|^2 = \frac{1}{n}\,\Var(e^{itZ_1}) =
\frac{1}{n}\,\big(1 - |\E\,e^{itZ_1}|^2\big) =
\frac{1}{n}\,\big(1 - e^{-t^2}\big).
$$
Hence, applying Plancherel's theorem and using the identity (4.7) for the
functions $\psi_r(\alpha)$ with $r = \alpha = 0$, we also have
\bee
\E\, \omega^2(G_n,\Phi) 
  & = &
\frac{1}{2\pi}\, \int_{-\infty}^\infty \E \, 
\Big|\frac{g_n(t) - g(t)}{t}\Big|^2\,dt \\
  & = &
\frac{1}{2\pi n} \int_{-\infty}^\infty 
\frac{1 - e^{-t^2}}{t^2}\,dt \, = \, \frac{1}{n \sqrt{\pi}}.
\ene
Thus, on average the $L^2$-distance 
$\omega(G_n,\Phi)$ is of order $1/\sqrt{n}$ as well.

Similar properties may be expected for the random variables 
$Z_k = \sqrt{n}\,\theta_k$ and hence for the random vector $X$.
Note that $|X| = \sqrt{n}$, while
$$
\E \left<X,\theta\right>^2 = \frac{1}{n}\,\sum_{k=1}^n \,
(\sqrt{n}\,\theta_k)^2 = 1, \quad \theta \in \S^{n-1},
$$
so that $X$ is isotropic. We now involve an asymptotic formula of 
Corollary 5.1 which yields
$$
\E_\theta\, \omega^2(F_\theta,\Phi)  \, = \,
\frac{1}{\sqrt{\pi}}\, \Big(1 + \frac{1}{4n}\Big)\,
\E\,\Big(1 - (1-\xi)^{1/2}\Big) - \frac{1}{8n\sqrt{\pi}} + 
O\Big(\frac{1}{n^2}\Big),
$$
where $\xi = \frac{\left<X,Y\right>}{n}$ with $Y$ being 
an independent copy of $X$. By the definition, $\xi$ takes only 
two values, 1 with probability $\frac{1}{n}$ and 0 with 
probability $1 - \frac{1}{n}$. Hence, the last expectation 
is equal to $\frac{1}{n}$, and we get
$$
\E_\theta\, \omega^2(F_\theta,\Phi)  \, = \,
\frac{7/8}{n \sqrt{\pi}} + O\Big(\frac{1}{n^2}\Big).
$$

As for the Kolmogorov distance, one may apply again Theorem 1.3,
which leads to the two-sided bound (13.3). Apparently, both 
logarithmic terms can be removed. Their appearance 
here is explained by the use of the Fourier tools (in the form 
of the Berry-Esseen bounds), while the proof of the 
Dvoretzky-Kiefer-Wolfowitz inequality on $\rho(G_n,\Phi)$ in
\cite{D-K-W} is based on the entirely different arguments.


\vskip5mm
\section{{\bf Improved Rates for Lacunary Systems}}
\setcounter{equation}{0}

\vskip2mm
\noindent
An orthonormal sequence of random variables 
$\{X_k\}_{k = 1}^\infty$ in $L^2(\Omega,{\mathfrak F},\P)$ 
is called a lacunary system of order $p>2$, if for any sequence
$(a_k)$ in $\ell^2$, the series $\sum_{k=1}^\infty a_k X_k$
converges in $L^p$-norm to an element of $L^p(\Omega,{\mathfrak F},\P)$. 
This property is equivalent to the validity of the Khinchine-type inequality
\be
\big(\E\,|a_1 X_1 + \dots + a_n X_n|^p\big)^{1/p} \, \leq \, M_p\,
(a_1^2 + \dots + a_n^2)^{1/2}
\en
for arbitrary $a_k \in \R$ with some constant $M_p$ independent of $n$ and 
the choice of the coefficients $a_k$. For basic properties of such systems
we refer an interested reader to the books \cite{K-Sa, Ka-S}.

Starting from an orthonormal lacunary system of order 
$p = 4$, consider the random vector $X = (X_1,\dots,X_n)$. 
According to Theorem 1.1, if $|X|^2 = n$ a.s. and $\E X = 0$, 
then
\be
c\,\E_\theta\, \omega^2(F_\theta,\Phi)  \, \leq \,
\frac{1}{n^3}\,\E \left<X,Y\right>^3 + 
\frac{1}{n^4}\,\E \left<X,Y\right>^4,
\en
where $Y$ is an independent copy of $X$. A similar bound
\be
c\,\E_\theta\,\rho^2(F_\theta,\Phi) \, \leq \,
\frac{\log n}{n^3}\, \E \left<X,Y\right>^3 + 
\frac{(\log n)^2}{n^4}\,\E \left<X,Y\right>^4
\en
also holds for the Kolmogorov distance. As easily follows from (15.1),
$$
\E\, |\left<X,Y\right>|^p \leq M_p^{2p} n^{p/2}.
$$ 
In particular,
$$
\E\, |\left<X,Y\right>|^3 \, \leq \, M_3^6\, n^{3/2}, \quad
\E \left<X,Y\right>^4 \, \leq \, M_4^8\, n^2.
$$
Hence, the bounds (15.2)-(15.3) lead to the estimates
\bee
c\,\E_\theta\, \omega^2(F_\theta,\Phi) 
 & \leq &
\frac{1}{n^{3/2}}\,M_3^6 + \frac{1}{n^2}\,M_4^8, \\
c\,\E_\theta\, \rho^2(F_\theta,\Phi) 
 & \leq &
\frac{\log n}{n^{3/2}}\,M_3^6 + 
\frac{(\log n)^2}{n^2}\,M_4^8.
\ene
Thus, if $M_4$ is bounded, both distances are at most of order 
$n^{-3/4}$ on average (modulo a logarithmic factor). Moreover, if
\be
\Sigma_3(n) \ \equiv \
\E \left<X,Y\right>^3 \ = \, \sum_{1 \leq i_1, i_2, i_3 \leq n}
\big(\E X_{i_1} X_{i_2} X_{i_3}\big)^2
\en
is bounded by a multiple of $n$, then these distances are 
on average at most $1/n$ (modulo a logarithmic factor
in the case of $\rho$). 

For an illustration, on the interval $\Omega = (-\pi,\pi)$ with 
the uniform measure $d\P(t) = \frac{1}{2\pi}\,dt$, consider 
a finite trigonometric system $X = (X_1,\dots,X_n)$ with 
components
\bee
X_{2k-1}(t) 
 & = &
\sqrt{2}\,\cos(m_k t), \\
X_{2k}(t) 
 & = &
\sqrt{2}\,\sin(m_k t), \qquad k = 1,\dots,n/2,
\ene
where $m_k$ are positive integers such that
$\frac{m_{k+1}}{m_k} \geq q > 1$ (assuming that $n$ is even).
Then $X$ is an isotropic random vector satisfying $|X|^2 = n$ 
and $\E X = 0$, and with $M_4$ bounded by a function of $q$ only.
For evaluation of the moment $\Sigma_3(n)$, one may use the identities
$$
\cos t = \E_\ep\, e^{i\ep t}, \quad 
\sin t = \frac{1}{i}\,\E_\ep\, \ep\,e^{i\ep t},
$$
where $\ep$ is a Bernoulli random variable taking the values $\pm 1$ 
with probability $\frac{1}{2}$. Let $\ep_1,\ep_2,\ep_3$ be independent 
copies of $\ep$. Using the property that $\ep_1 \ep_3$ and $\ep_2 \ep_3$ 
are independent, the first identity implies that, for all integers 
$1 \leq n_1 \leq n_2 \leq n_3$,
\bee
\E\,\cos(n_1 t) \cos(n_2 t) \cos(n_3 t) 
 & = &
\E_\ep\,\E\,\exp\{i(\ep_1 n_1 + \ep_2 n_2 + \ep_3 n_3)\,t\} \\
 & = &
\E_\ep\,I\{\ep_1 n_1 + \ep_2 n_2 + \ep_3 n_3 = 0\} \\
 & = &
\E_\ep\,I\{\ep_1 n_1 + \ep_2 n_2 = n_3\} \ = \  
\frac{1}{4}\,I\{n_1 + n_2 = n_3\},
\ene
where $\E_\ep$ means the expectation over $(\ep_1,\ep_2,\ep_3)$, 
and where $I\{A\}$ denotes the indicator of the event $A$. Similarly, 
involving also the identity for the sine function, we have
\bee
\E\,\sin(n_1 t) \sin(n_2 t) \cos(n_3 t) 
 & = &
-\E_\ep\,\E\,\ep_1 
\ep_2\, \exp\{i(\ep_1 n_1 + \ep_2 n_2 + \ep_3 n_3)\,t\} \\
 & & \hskip-20mm = \
-\E_\ep\,
\ep_1 \ep_2\,I\{\ep_1 n_1 + \ep_2 n_2 + \ep_3 n_3 = 0\} \\
 & & \hskip-20mm = \
-\E_\ep\,
\ep_1 \ep_2\,I\{\ep_1 n_1 + \ep_2 n_2 = n_3\}
 \ = \ 
-\frac{1}{4}\,I\{n_1 + n_2 = n_3\},
\ene
\bee
\E\,\sin(n_1 t) \cos(n_2 t) \sin(n_3 t) 
 & = &
-\E_\ep\,\E\,
\ep_1 \ep_3\, \exp\{i(\ep_1 n_1 + \ep_2 n_2 + \ep_3 n_3)\,t\} \\
 & & \hskip-20mm = \
-\E_\ep\,
\ep_1 \ep_3\,I\{\ep_1 n_1 + \ep_2 n_2 + \ep_3 n_3 = 0\} \\
 & & \hskip-20mm = \
-\E_\ep\,
\ep_1\,I\{\ep_1 n_1 + \ep_2 n_2 = n_3\} \ = \  
-\frac{1}{4}\,I\{n_1 + n_2 = n_3\},
\ene
\bee
\E\,\cos(n_1 t) \sin(n_2 t) \sin(n_3 t) 
 & = &
-\E_\ep\,\E\,
\ep_2 \ep_3\, \exp\{i(\ep_1 n_1 + \ep_2 n_2 + \ep_3 n_3)\,t\} \\
 & & \hskip-20mm = \
-\E_\ep\,
\ep_2 \ep_3\,I\{\ep_1 n_1 + \ep_2 n_2 + \ep_3 n_3 = 0\} \\
 & & \hskip-20mm = \
-\E_\ep\,
\ep_2\,I\{\ep_1 n_1 + \ep_2 n_2 = n_3\} \ = \ 
- \frac{1}{4}\,I\{n_1 + n_2 = n_3\}.
\ene
On the other hand, if the sine function appears in the product once
or three times, such expectations will be vanishing. They are thus 
vanishing in all cases where $n_1 + n_2 \neq n_3$, and do not 
exceed $\frac{1}{4}$ in absolute value for any combination of 
sine and cosine terms in all cases with $n_1 + n_2 = n_3$. 
Therefore, the moment $\Sigma_3(n)$ in (15.4) is bounded by
a multiple of
$$
T_3(n) \, = \,
{\rm card}\big\{(i_1,i_2,i_3): 1 \leq i_1 \leq i_2 < i_3 \leq n, \ 
m_{i_1} + m_{i_2} = m_{i_3}\big\}.
$$

One can now involve the lacunary assumption. If $q \geq 2$,
the property $i_1 \leq i_2 < i_3$ implies 
$m_{i_1} + m_{i_2} < m_{i_3}$, so that $T_3(n) = \Sigma_3(n) = 0$.
In the case $1 < q < 2$, define $A_q$ to be the (finite) collection 
of all couples $(k_1,k_2)$ of positive integers such that 
$$
q^{-k_1} + q^{-k_2} \geq 1.
$$
By the lacunary assumption, if $1 \leq i_1 \leq i_2 < i_3 \leq n$,
we have 
$$
m_{i_1} + m_{i_2} \, \leq \, \big(q^{-(i_3 - i_1)} + 
q^{-(i_3 - i_2)}\big)\,m_{i_3} \, < \, m_{i_3},
$$
as long as the couple $(i_3 - i_1,i_2 - i_1)$ is not in $A_q$.
Hence, 
\bee
T_3(n) 
 & \leq &
{\rm card}\big\{(i_1,i_2,i_3): 1 \leq i_1 \leq i_2 < i_3 \leq n, \ 
(i_3 - i_1,i_2 - i_1) \in A_q\big\} \\
 & \leq &
n\,{\rm card}(A_q) \ \leq \ c_q n
\ene
with constant depending on $q$ only. Returning to (15.2)-(15.3), 
we then obtain:

\vskip5mm
{\bf Proposition 15.1.} {\sl For the lacunary trigonometric 
system $X$ of an even length $n$ and with parameter $q>1$, we have
$$
\E_\theta\, \omega^2(F_\theta,\Phi)  \, \leq \, \frac{c_q}{n^2}, \quad
\E_\theta\,\rho^2(F_\theta,\Phi) \, \leq \, \frac{c_q\,(\log n)^2}{n^2},
$$
where the constants $c_q$ depend $q$ only. 
}

\vskip5mm
In this connection one should mention a classical result of Salem
and Zygmund concerning distributions of the lacunary sums
$$
S_n = \sum_{k=1}^n \, (a_k \cos(m_k t) + b_k \sin(m_kt))
$$
with an arbitrary prescribed sequence of the coefficients 
$(a_k)_{k \geq 1}$ and $(b_k)_{k \geq 1}$. Assume that 
$\frac{m_{k+1}}{m_k} \geq q > 1$ for all $k$ and put 
$$
v_n^2 \, = \, \frac{1}{2}\, \sum_{k=1}^n \, (a_k^2 + b_k^2) 
\qquad (v_n \geq 0),
$$
so that the normalized sums $Z_n = S_n/v_n$ have mean zero 
and variance one under the measure $\P$. It was shown in 
\cite{S-Z1} that $Z_n$ are weakly convergent to the standard 
normal law, i.e., their distributions $F_n$ under $\P$ satisfy
$\rho(F_n,\Phi) \rightarrow 0$ as $n \rightarrow \infty$,
if and only if $\frac{a_n^2 + b_n^2}{v_n^2} \rightarrow 0$
(in fact, the weak convergence was established on every subset 
of $\Omega$ of positive measure).

Restricting to the coefficients $\theta_{2k-1} = a_k/v_n$, 
$\theta_{2k} = b_k/v_n$, Salem-Zygmund's theorem may be
stated as the assertion that $\rho(F_\theta,\Phi)$ is small, 
if and only if $\|\theta\|_\infty = \max_{1 \leq k \leq n} |\theta_k|$ 
is small. The latter condition naturally appears in the central limit 
theorem for weighted sums of independent identically distributed 
random variables. Thus, Proposition 15.1 complements this result 
in terms of the rate of convergence in the mean on the unit sphere.
It would be interesting to describe explicit coefficients $\theta_k$,
for which we get a standard rate of normal approximation
(perhaps, using other approaches such as the Stein method,
cf. e.g. \cite{G-R}).

The result of \cite{S-Z1} was generalized in \cite{S-Z2};
it turns out there is no need to assume that all $m_k$ are integers, 
and the asymptotic normality is preserved for real $m_k$ such that 
$\inf_k \frac{m_{k+1}}{m_k} > 1$. However, in this more 
general situation, the rate $1/n$ as in Proposition 15.1 is 
no longer true (although the rate $1/\sqrt{n}$ is valid). 
The main reason is that the means
$$
\E X_{2k-1} \, = \, \sqrt{2}\,\E\,\cos(m_k t) \, = \, 
\sqrt{2}\ \frac{\sin(\pi m_k)}{\pi m_k}
$$
may be non-zero. For example, choosing 
$m_k = 2^k + \frac{1}{2}$, we obtain an orthonormal system
with $\E X_{2k} = 0$, while 
$$
\E X_{2k-1} = \frac{2\sqrt{2}}{\pi\, (2^{k+1} + 1)}.
$$
Hence
$$
\E \left<X,Y\right> \, = \, |\E X|^2 = \frac{8}{\pi^2}\, 
\sum_{k=1}^n \frac{1}{(2^{k+1} + 1)^2} \, \rightarrow \, c 
\qquad (n \rightarrow \infty)
$$
for some absolute constant $c>0$ (where $Y$ is 
an independent copy of $X$). In this situation, as was
already mentioned in (5.3), cf. Remark 5.3, we have a lower bound
$$
\E_\theta\, \omega^2(F_\theta,F) \geq 
\frac{c}{2\sqrt{\pi}\,n} + O\Big(\frac{1}{n^2}\Big).
$$
Since $\E \left<X,Y\right>^3 = O(n)$ and 
$\E \left<X,Y\right>^4 = O(n^2)$, this inequality may actually 
be replaced with equality, according to (5.2).
A similar asymptotic holds as well when $F$ is replaced with~$\Phi$.


\vskip5mm
\section{{\bf Improved Rates for Independent and Log-concave Summands}}
\setcounter{equation}{0}

\vskip2mm
\noindent
Let $X = (X_1,\dots,X_n)$ be an isotropic random vector in $\R^n$
with mean zero. If the components $X_k$ are independent, the normal 
approximation for the distributions $F_\theta$ of the weighted sums 
$$
S_\theta = \theta_1 X_1 + \dots + \theta_n X_n, \quad \theta \in \S^{n-1},
$$
may be controlled by virtue of the Berry-Esseen theorem under the
3-rd moment assumption. Namely, this theorem provides 
an upper bound
\be
\rho(F_\theta,\Phi) \leq c \sum_{i=1}^n |\theta_i|^3\,\E\,|X_i|^3
\en
(cf. e.g. \cite{P1}, \cite{P2}). Since $\E\,|X_i|^3 \geq 1$, 
the sum in (16.1) is at least $\frac{1}{\sqrt{n}}$. On the other hand, 
(16.1) yields an upper estimate on average
\be
\E_\theta\, \rho(F_\theta,\Phi) \leq \frac{c\beta_3}{\sqrt{n}}, \quad
\beta_3 = \max_{1 \leq i \leq n} \E\,|X_i|^3,
\en
which is consistent with the standard rate.

As it turns out, the relations (16.1)-(16.2) are far from being optimal
for most of $\theta$, as the following statement due to Klartag and Sodin shows. 

\vskip5mm
{\bf Theorem 16.1} (\cite{K-S}). {\sl If the random variables $X_1,\dots,X_n$
are independent, have mean zero, variance one, and finite 4-th moments,
then
\be
\E_\theta\, \rho(F_\theta,\Phi) \leq \frac{c\beta_4}{n}, \quad 
\beta_4 = \frac{1}{n} \sum_{i=1}^n \E X_i^4.
\en
Moreover, for any $r \geq 0$,
$$
\s_{n-1}\big\{n\rho(F_\theta,\Phi) \geq c\beta_4 r\big\} \, \leq \, 2\, e^{-\sqrt{r}}.
$$
}

In the i.i.d. case, $\beta_4 = \E X_1^4$, and we obtain 
an upper bound of order at most $1/n$.

In fact, in the i.i.d. case, the relation (16.3) may be further
sharpened under the 5-th moment assumption, if $\E X_1^3 = 0$,
and if $\Phi(x)$ is slightly modified to
$$
G(x) = \Phi(x) - \frac{\beta_4 - 3}{8(n+2)}\,(x^3 - 3x)\,\varphi(x),
\quad x \in \R,
$$
where $\varphi(x) = \frac{1}{\sqrt{2\pi}}\,e^{-x^2/2}$ is the
standard normal density.

\vskip5mm
{\bf Theorem 16.2.} {\sl If the random variables $X_1,\dots,X_n$
are independent, identically distributed, and have moments 
$\E X_1 = 0$, $\E X_1^2 = 1$, $\E X_1^3 = 0$, $\E X_1^4 = \beta_4$,
$\E\,|X_1|^5 = \beta_5 < \infty$, then
\be
\E_\theta\, \rho(F_\theta,G) \leq \frac{c\beta_5}{n^{3/2}}.
\en
Moreover, for any $r \geq 0$,
$$
\s_{n-1}\Big\{n^{3/2}\rho(F_\theta,G) \geq c\beta_4 r\Big\} \, \leq \, 
2\, \exp\{-r^{2/5}\}.
$$
}

\vskip2mm
We refer an interested reader to \cite{B3} and \cite{B-C-G6}. In the i.i.d. 
case, both inequalities (16.3) and (16.4) are sharp in the following sense. 
If $\alpha_3 = \E X_1^3 \neq 0$ and $\beta_4 < \infty$, then, for any 
function $G$ of bounded total variation, such that $G(-\infty) = 0$ and
$G(\infty) = 1$, we have
$$
\E_\theta\, \rho(F_\theta,G) \geq \frac{c}{n}
$$
with a constant $c>0$ depending on $\alpha_3$ and $\beta_4$.
Similarly, if $\alpha_3 = 0$, $\beta_4 \neq 3$, $\beta_5 < \infty$, then
$$
\E_\theta\, \rho(F_\theta,G) \geq \frac{c}{n^{3/2}},
$$
where the constant $c>0$ depends on $\beta_4$ and $\beta_5$ only.

In the upper bounds such as (16.3), the independence assumption may be replaced
with closely related hypotheses. The random vector $X$ is said to have 
a log-concave distribution, when it has a density of the form $p(x) = e^{-V(x)}$ 
where $V:\R^n \rightarrow (-\infty,\infty]$ is a convex function. Recall that
the distribution of $X$ is coordinatewise symmetric, if
$$
p(\ep_1 x_1,\dots,\ep_n x_n) = p(x_1,\dots,x_n), \quad x_i \in \R,
$$
for any choice of signs $\ep_i = \pm 1$. The following theorem sharpening 
(16.1) is due to Klartag.

\vskip5mm
{\bf Theorem 16.3} (\cite{Kl}). {\sl Suppose that the isotropic random vector 
$X = (X_1,\dots,X_n)$ in $\R^n$ has a coordinatewise symmetric log-concave 
distribution. For all $\theta = (\theta_1,\dots,\theta_n) \in \S^{n-1}$,
\be
\|F_\theta - \Phi\|_{\rm TV} \leq c \sum_{i=1}^n \theta_i^4.
\en
}

\vskip2mm
Here, the total variation distance is understood in the usual sense as
$$
\|F_\theta - \Phi\|_{\rm TV} = \int_{-\infty}^\infty
|p_\theta(x) - \varphi(x)|\,dx,
$$
where $p_\theta$ denotes the density of $S_\theta$. By the assumptions, 
$p_\theta$ is symmetric about the origin and is log-concave for any 
$\theta \in \S^{n-1}$. Note that, by the coordinatewise symmetry, 
the isotropy assumption is reduced to the moment condition 
$\E X_i^2 = 1$ ($1 \leq i \leq n$).

In particular, it follows from (16.5) that
\be
\E_\theta\, \rho(F_\theta,\Phi) \leq 
\E_\theta\, \|F_\theta - \Phi\|_{\rm TV} \leq \frac{c}{n}.
\en


\vskip5mm
\section{{\bf Improved Rates Under Correlation-Type Conditions}}
\setcounter{equation}{0}

\vskip2mm
\noindent
Up to a logarithmically growing term, the improved rate as in the upper
bound (16.3) can be achieved under more flexible correlation-type conditions
(in comparison with independence). For example, one may consider
an optimal value $\Lambda = \Lambda(X)$ in the relation
\be
\Var\bigg(\sum_{i,j=1}^n a_{ij} X_i X_j\bigg) \leq 
\Lambda \sum_{i,j=1}^n a_{ij}^2 \qquad (a_{ij} \in \R),
\en
which we call that the random vector $X = (X_1,\dots,X_n)$ 
satisfies a second order correlation condition with constant $\Lambda$. 
This quantity is finite as long as the moment $\E\,|X|^4$ is finite. 

To relate $\Lambda$ to the moment-type characteristics which we discussed
before, one may apply (17.1) with $a_{ij} = \delta_{ij}$ or (as another option)
with $a_{ij} = \theta_i \theta_j$,
$\theta = (\theta_1,\dots,\theta_n) \in \S^{n-1}$. This gives that
$$
\sigma_4^2 \leq \Lambda, \quad m_4^2 \leq 
\sup_{\theta \in \S^{n-1}} \E S_\theta^4 \leq 1 + \Lambda,
$$
where in the last inequality we should assume that 
$\E S_\theta^2 = 1$ for all $\theta$ (i.e. $X$ is isotropic).
In the latter case, necessarily $\Lambda \geq \frac{n-1}{n}$, so that
$\Lambda$ is bounded away from zero.

If the distribution of $X$ is ``regular" in some sense, one may also
bound $\Lambda$ from above. For example, this is the case when
it shares a Poincar\'e-type inequality
\be
\lambda_1 \Var(u(X)) \leq \E\,|\nabla u(X)|^2,
\en
which is required to hold in the class of all bounded, smooth functions 
$u$ on $\R^n$ with a constant $\lambda_1 > 0$ independent of $u$
(called the spectral gap). We then have
\be
\Lambda \leq \frac{4}{\lambda_1^2}, \quad \Lambda \leq \frac{4}{\lambda_1},
\en
where in the second inequality we assume that $X$ is isotropic.

The following relation is established in \cite{B-C-G4}.

\vskip5mm
{\bf Theorem 17.1.} {\sl If the distribution of $X$ is 
isotropic and symmetric about the origin, then
\be
\E_\theta\,\rho(F_\theta,\Phi) \leq c\Lambda\,\frac{\log n}{n}.
\en
}

\vskip2mm
The proof is based on the second order spherical concentration phenomenon 
which was developed in \cite{B-C-G1} with the aim of applications to 
randomized central limit theorems. It indicates that the deviations of 
any smooth function $u(\theta)$ on $\S^{n-1}$ from the mean 
$\E_\theta u(\theta)$ are at most of the order $1/n$, provided that 
$u$ is orthogonal in $L^2(\R^n,\s_{n-1})$ to all linear functions and 
has a ``bounded" Hessian (the matrix of second order partial derivatives). 
Being applied to the characteristic functions $u(\theta) = f_\theta(t)$, 
this property yields an upper bound
$$
\E_\theta\, |f_\theta(t) - f(t)|^2\,\leq\,\frac{c \Lambda t^4}{n^2}
$$
on every interval $|t| \leq An^{1/5}$ with constants $c > 0$ 
depending on the parameter $A \geq 1$ only. This estimate
can be used to bound the integrals in (8.4) to get a similar 
variant of (17.4).

The symmetry hypothesis in Proposition 17.1 may be dropped, if $\Lambda$
is replaced by $\lambda_1^{-1}$ which is a larger quantity according
to (17.3). In addition, one can control large deviations of the distance 
$\rho(F_\theta,\Phi)$ for most of the directions $\theta$ (rather than on average).
The corresponding assertions are obtained in \cite{B-C-G5}.

\vskip5mm
{\bf Theorem 17.2.} {\sl Let $X$ be an isotropic random vector 
in $\R^n$ with mean zero and a positive Poincar\'e constant 
$\lambda_1$. Then
\be
\E_\theta\, \rho(F_\theta,\Phi) \, \leq \, 
c\lambda_1^{-1}\,\frac{\log n}{n}.
\en
Moreover, for all $r > 0$,
$$
\s_{n-1}\Big\{\rho(F_\theta,\Phi) \geq c\lambda_1^{-1} 
\frac{\log n}{n}\, r\Big\} \, \leq \, 2\,e^{-\sqrt{r}}.
$$
}

\vskip2mm
The logarithmic term in (17.5) may be removed using the
less sensitive $L^2$-distance:
$$
\E_\theta\, \omega^2(F_\theta,\Phi) \leq \frac{c}{\lambda_1^2\, n^2}.
$$

There is an extensive literature devoted to bounding the
spectral gap $\lambda_1$ from below. In particular, it is positive
for any log-concave probability distribution on $\R^n$.
A well-known conjecture raised by Kannan, Lov\'asz and Simonovits
asserts that $\lambda_1$ is actually bounded away from zero, as long as
the random vector $X$ has an isotropic log-concave distribution
(cf. \cite{K-L-S}). The best known dimensional lower bound up to date
is due to Klartag and Lehec \cite{K-L} who showed that
$$
\lambda_1 \geq \frac{c}{(\log n)^\alpha}
$$
for some absolute positive constants $c$ and $\alpha$ (one may take
$\alpha = 10$). Applying this bound in Theorem 17.2, we therefore obtain:

\vskip5mm
{\bf Corollary 17.3.} {\sl Let $X$ be an isotropic random vector 
in $\R^n$ with mean zero and a log-concave probability distribution. Then
with some absolute positive constants $c$ and $\alpha$
\be
\E_\theta\, \rho(F_\theta,\Phi) \, \leq \, 
\frac{c(\log n)^\alpha}{n}.
\en
}

\vskip2mm
Thus, there is a certain extension of Klartag's bound (16.6)
at the expense of a logarithmic factor to the entire
class of isotropic log-concave probability distributions on $\R^n$.

One may also argue in the opposite direction:
upper bounds of the form 
$$
\E_\theta\, \rho(F_\theta,\Phi) \leq \frac{c(\log n)^\beta}{n}, \quad
\beta > 0,
$$
in the class of log-concave probability distributions on $\R^n$
imply lower bounds $\lambda_1 \geq c\,(\log n)^{-\beta'}$ with some 
$\beta'>0$, cf. \cite{B-C-G4}.

\vskip 0.5cm

\address
{Sergey G. Bobkov,}
\email {bobkov@math.umn.edu}
 
\address
{Gennadiy P. Chistyakov,
\email {chistyak@math.uni-bielefeld.de}

\address
{Friedrich G\"otze,
}
\email {goetze@math.uni-bielefeld.de}

\end{document}

\bibitem{C}
Y. Chen. An almost constant lower bound of the isoperimetric coefficient 
      in the KLS conjecture. Geom. Funct. Anal. 31 (2021), no. 1, 34--61.

\vskip 0.5cm

\address
{Sergey G. Bobkov \newline
School of Mathematics, University of Minnesota  \newline 
127 Vincent Hall, 206 Church St. S.E., Minneapolis, MN 55455 USA
\smallskip}
\email {bobkov@math.umn.edu}

\address
{Gennadiy P. Chistyakov\newline
Fakult\"at f\"ur Mathematik, Universit\"at Bielefeld\newline
Postfach 100131, 33501 Bielefeld, Germany}
\email {chistyak@math.uni-bielefeld.de}

\address
{Friedrich G\"otze\newline
Fakult\"at f\"ur Mathematik, Universit\"at Bielefeld\newline
Postfach 100131, 33501 Bielefeld, Germany}
\email {goetze@mathematik.uni-bielefeld.de}

\bibitem{B-C-G6}
S. G. Bobkov, G. P. Chistyakov, and F. G\"otze.  Concentration and Gaussian 
      approximation for randomized sums. In preparation.